\documentclass[twoside]{amsart}
\usepackage[centertags]{amsmath}
\usepackage{amsfonts}
\usepackage{amssymb}
\usepackage{amsthm}
\usepackage[ansinew]{inputenc}
\usepackage[cmtip,all]{xy}
\usepackage{graphicx}
\newtheorem{thm}[equation]{Theorem}
\newtheorem{cor}[equation]{Corollary}
\newtheorem{lem}[equation]{Lemma}
\newtheorem{prop}[equation]{Proposition}
\theoremstyle{definition}
\newtheorem{defn}[equation]{Definition}
\theoremstyle{remark}
\newtheorem{rem}[equation]{Remark}

\numberwithin{equation}{section}
\newcommand{\norm}[1]{\left\Vert#1\right\Vert}

\newcommand{\set}[1]{\left\{#1\right\}}
\newcommand{\Real}{\mathbb R}

\newcommand{\To}{\longrightarrow}

\newcommand{\C}[1]{\mathbf{#1}} 

\def\ni{{nil}}
\def\abb{{ab}}

\def\r{\rightarrow} 
\def\rr{\Rightarrow} 

\def\id{\operatorname{Id}}
\newcommand{\sym}[1]{\operatorname{Sym}(#1)}

\newcommand{\symt}[1]{\operatorname{Sym}_\vc(#1)}
\newcommand{\symtex}[1]{\overline{\operatorname{Sym}}_\vc(#1)}

\def\sign{\operatorname{sign}}

\def\aut{\operatorname{Aut}}

\def\endo{\operatorname{End}}
\def\enda{\operatorname{\mathbf{End}}}

\def\st{\stackrel} 

\def\ul{\underline} 

\newcommand{\hopf}{\mathit{Hopf}}

\newcommand{\Ab}{\mathrm{\mathbf{Ab}}}

\def\Z{\mathbb{Z}}

\def\S{\Sigma}

\def\L{\Omega}

\newcommand{\grupo}[1]{\langle #1\rangle}
\newcommand{\grupoo}[1]{\langle\!\langle #1\rangle\!\rangle}

\newcommand{\ad}{\mathsf{Ad}}

\newcommand{\vc}{\Box}    
\newcommand{\vi}{\boxminus}  

\begin{document}

\title[The symmetric action on secondary homotopy groups]{The symmetric action on secondary homotopy groups}%
\author{Hans-Joachim Baues and Fernando Muro}%
\address{Max-Planck-Institut f\"ur Mathematik, Vivatsgasse 7, 53111 Bonn, Germany}%
\email{baues@mpim-bonn.mpg.de, muro@mpim-bonn.mpg.de}%

\thanks{The second author was partially supported
by the project MTM2004-01865 and the MEC postdoctoral fellowship EX2004-0616.}%

\begin{abstract}
We show that the symmetric track group $\symt{n}$, which is an extension of the symmetric group $\sym{n}$
associated to the second Stiefel-Withney class, acts as a crossed
module on the secondary homotopy group of a pointed
space.
\end{abstract}
\maketitle

\section*{Introduction}

Secondary homotopy operations like Toda brackets \cite{toda} or cup-one-products \cite{K1p}, \cite{c1Tb}, are defined by pasting
tracks, where tracks are homotopy classes of homotopies. Since secondary homotopy operations play a crucial role in homotopy
theory it is of importance to develop the algebraic theory of tracks. We do this by introducing secondary
homotopy groups of a pointed space $X$
$$\Pi_{n,*}X=\left(\Pi_{n,1}X\st{\partial}\r\Pi_{n,0}X\right)$$
which have the structure of a quadratic pair module, see Section \ref{sqpm}. 
Here $\partial$ is a group homomorphism
with cokernel $\pi_nX$ and kernel $\pi_{n+1}X$ for $n\geq 3$. 

We define $\Pi_{n,*}X$ for $n\geq 2$ directly in terms of maps 
$S^n\r X$ and tracks from such maps to the trivial map. For $n\geq 0$ the functor $\Pi_{n,*}$ is an additive
version of the functor $\pi_{n,*}$ studied in \cite{2hg1}. The homotopy category of $(n-1)$-connected
$(n+1)$-types is equivalent via $\Pi_{n,*}$ to the homotopy category of quadratic pair modules for $n\geq 3$.

In this paper we consider the ``generalized coefficients'' of secondary homotopy groups $\Pi_{n,*}X$ obtained 
by the action of the symmetric group $\sym{n}$ on $S^n=S^1\wedge\st{n}\cdots\wedge
S^1$ via permutation of coordinates. For a permutation $\sigma\in\sym{n}$ the map $\sigma\colon S^n\r S^n$ has
degree $\sign \sigma\in\set{\pm1}$. The group $\set{\pm1}$ also acts on $S^n$ by using the
topological abelian group structure of $S^1$ and suspending $(n-1)$ times. This shows that there are tracks
$\sigma\rr\sign\sigma$ which by definition are the elements of the symmetric track group $\symt{n}$. Also these tracks act on
$\Pi_{n,*}X$.
We clarify this action by showing that the group $\symt{n}$ gives rise to a
crossed module which acts as a crossed module on the quadratic pair module $\Pi_{n,*}X$.



The symmetric track group is a central extension
$$\Z/2\hookrightarrow\symt{n}\st{\delta}\twoheadrightarrow\sym{n}$$
which, as we show, represents the second Stiefel-Withney class pulled back to $\sym{n}$. 
The symmetric track group is computed in Section \ref{krio}. We actually compute a
faithful positive pin representation of $\symt{n}$ from which we derive a finite presentation of this group. 
This group also
arose in a different way in the work of Schur \cite{dsag} and Serre \cite{iwf}.

In \cite{2hg3} we describe the smash product operation on 
secondary homotopy groups $\Pi_{n,*}X$. This operation 
endows $\Pi_{*,*}$ with the structure of a lax symmetric monoidal functor where the crossed module action of
$\symt{n}$ on $\Pi_{n,*}X$ is of crucial importance. 
This leads to an algebraic approximation of the symmetric monoidal 
category of spectra by secondary homotopy groups, see \cite{2hg4}.
As an example we prove a formula for the unstable cup-one-product
$\alpha\cup_1\alpha\in\pi_{2n+1}S^{2m}$ of an element $\alpha\in\pi_nS^m$ where $n$ and $m$ are even. We show
that $$2(\alpha\cup_1\alpha)=\frac{n+m}{2}(\alpha\wedge\alpha)(\S^{2(n-1)}\eta)$$
where $\eta\colon S^3\r S^2$ is the Hopf map. If $n/2$ is odd and $m/2$ is even then this formula was achieved by
totally different methods in \cite{K1p}.

\section{Square groups and quadratic pair modules}\label{sqpm}

In this section we describe the algebraic concepts needed for the structure of secondary homotopy groups. 

\begin{defn}
A \emph{square group} $X$ is a diagram
$$X=(X_e\mathop{\leftrightarrows}\limits^P_HX_{ee})$$
where $X_e$ is a group with an additively written group law, $X_{ee}$
is an abelian group, $P$ is a homomorphism, $H$ is a function such
that the crossed effect
$$(a|b)_H=H(a+b)-H(b)-H(a)$$
is linear in $a$ and $b\in X_e$, and the following relations are
satisfied, $x,y\in X_{ee}$,
\begin{enumerate}
\item $(Px|b)_H=0$, $(a|Py)=0$,
\item $P(a|b)_H=-a-b+a+b$,
\item $PHP(x)=P(x)+P(x)$.
\end{enumerate}

The function $$T=HP-1\colon X_{ee}\To X_{ee}$$
is an involution, i. e. a homomorphism with $T^2=1$.

A morphism of square groups $f\colon X\r Y$ is given by
homomorphisms
$$f_e\colon X_e\To Y_e,$$
$$f_{ee}\colon X_{ee}\To Y_{ee},$$
commuting with $P$ and $H$.

Let $\C{SG}$ be the category of square groups. A square group $X$
with $X_{ee}=0$ is the same as an abelian group $X_e$. This yields
the full inclusion of categories $\Ab\subset\C{SG}$ where $\Ab$ is the category of abelian groups.
\end{defn}

Square groups were introduced in \cite{ecg} to describe quadratic
endofunctors of the category $\C{Gr}$ of groups. More precisely, any
square group $X$ gives rise to a quadratic functor
$$-\otimes X\colon\C{Gr}\To\C{Gr}.$$
Given a group $G$ the group $G\otimes X$ is generated by the symbols
$g\otimes x$ and $[g,h]\otimes z$, $g,h\in G$, $x\in X_e$, $z\in
X_{ee}$ subject to the  relations
\begin{eqnarray*}
(g+h)\otimes x&=&g\otimes x+h\otimes x+[g,h]\otimes H(x),\\
{[}g,g{]}\otimes z&=&g\otimes P(z),
\end{eqnarray*}
where $g\otimes x$ is linear in $x$ and $[g,h]\otimes z$ is central and linear in each variable $g, h, z$.
If $X$ is an abelian group then $G\otimes X=G_\abb\otimes X_e$. In
fact, any quadratic functor $F\colon\C{Gr}\r\C{Gr}$ which preserves
reflexive coequalizers and filtered colimits has the form
$F=-\otimes X$, see \cite{ecg}.  The theory of square groups is discussed in detail in \cite{qaI}.

There is a natural isomorphism
$$X_e\st{\cong}\To \Z\otimes X,\;\; x\mapsto 1\otimes x.$$
In particular the homomorphism $n\colon \Z\r\Z$ induces a
homomorphism $n^*\colon X_e\r X_e$ fitting into the following
commutative diagram
\begin{equation}\label{n*}
\xymatrix{\Z\otimes X\ar[r]^{n\otimes
X}\ar@{<-}[d]_{\cong}&\Z\otimes
X\ar@{<-}[d]^{\cong}\\X_e\ar[r]^{n^*}&X_e}
\end{equation}
The homomorphism $n^*$ is explicitly given by the following formula,
\begin{equation*}
n^*x=n\cdot a+\binom{n}{2}PH(x). 
\end{equation*}
Here we set $\binom{n}{2}=\frac{n(n-1)}{2}$ and for any additively written group $G$ and any $n\in\Z$, $g\in G$, 
$$n\cdot g=\left\{
    \begin{array}{ll}
      g+\st{n}\cdots+g, & \hbox{if $n\geq0$;} \\
      &\\
      -g-\st{-n}\cdots-g, & \hbox{if $n<0$.}
    \end{array}
  \right.
$$


The function $n\cdot\,\colon G\r G$ in general is not a
homomorphism, but if $G$ is abelian then $n\cdot\,$ is a
homomorphism. This homomorphism is generalized by $n^*$ in (\ref{n*}) for square
groups.

\begin{defn}\label{qpm}
A \emph{quadratic pair module} $C$ is a morphism $\partial\colon
C_{(1)}\r C_{(0)}$ between square groups
\begin{eqnarray*}
C_{(0)}&=&(C_{0}\mathop{\leftrightarrows}^{P_0}_HC_{ee}),\\
C_{(1)}&=&(C_{1}\mathop{\leftrightarrows}^P_{H_1}C_{ee}),
\end{eqnarray*}
such that $\partial_{ee}=1\colon C_{ee}\r C_{ee}$ is the identity
homomorphism. In particular $\partial$ is completely determined by
the diagram
\begin{equation}\label{ya}
\xymatrix{&C_{ee}\ar[ld]_P&\\C_{1}\ar[rr]_\partial&&C_{0}\ar[lu]_H}
\end{equation}
where $\partial=\partial_e$, $H_1=H\partial$ and $P_0=\partial P$.

Morphisms of quadratic pair modules $f\colon C\r D$ are therefore
given by group homomorphisms
$f_{0}\colon C_{0}\r D_{0}$,
$f_{1}\colon C_{1}\r D_{1}$,
$f_{ee}\colon C_{ee}\r D_{ee}$,
commuting with $H$, $P$ and $\partial$ in (\ref{ya}) as in the diagram
$$\xymatrix{C_0\ar[r]^H\ar[d]^{f_0}&C_{ee}\ar[d]^{f_{ee}}\ar[r]^P&C_1\ar[d]^{f_1}\ar[r]^\partial&C_0\ar[d]^{f_0}\\
D_0\ar[r]^H&D_{ee}\ar[r]^P&D_1\ar[r]^\partial&D_0}$$
They form a
category denoted by $\C{qpm}$.

Quadratic pair modules are also the objects of a bigger category
$\C{wqpm}$ given by weak morphisms. A \emph{weak morphism} $f\colon
C\r D$ between quadratic pair modules is given by three
homomorphisms $f_0, f_1, f_{ee}$ as above, but we only require the
following two diagrams to be commutative
$$\begin{array}{cc}
\xymatrix{C_{ee}\ar[r]^T\ar[d]^{f_{ee}}&C_{ee}\ar[d]^{f_{ee}}\\D_{ee}\ar[r]^T&D_{ee}}
&
\xymatrix{\otimes^2(C_0)_\abb\ar[r]^<(.3){(-|-)_H}\ar[d]^{\otimes^2(f_0)_\abb}&C_{ee}\ar[d]^{f_{ee}}\ar[r]^P&C_1\ar[d]^{f_1}\ar[r]^\partial&C_0\ar[d]^{f_0}\\
\otimes^2(D_0)_\abb\ar[r]^<(.3){(-|-)_H}&D_{ee}\ar[r]^P&D_1\ar[r]^\partial&D_0}
\end{array}$$
Here $\otimes^2A=A\otimes A$ denotes the tensor square of an abelian group.
Therefore $\C{qpm}\subset\C{wqpm}$ is a subcategory with the same objects.

\end{defn}


Let $(\Z,\cdot)$ be the multiplicative (abelian) monoid of the
integers $\Z$.

\begin{defn}\label{sqmact}
Any quadratic pair module $C$ admits an action of $(\Z,\cdot)$ given
by the morphisms $n^*\colon C\r C$ in $\C{wqpm}$, $n\in\Z$,
defined by the equations
\begin{itemize}
\item $n^*x=n\cdot x+\binom{n}{2}\partial PH(x)$ for $x\in C_{0}$,
\item $n^*y=n\cdot y+\binom{n}{2}PH\partial(y)$ for $y\in C_{1}$,
\item $n^*z=n^2z$ for $z\in C_{ee}$.
\end{itemize}
We point out that $n^*\colon C\r C$ is an example of a weak morphism
which is not a morphism in $\C{qpm}$ since $n^*$ is not compatible
with $H$. Notice that $n^*\colon C_{0}\r C_{0}$ and $n^*\colon
C_{1}\r C_{1}$ are induced by   the
square group morphisms $n\otimes
C_{(0)}$ and $n\otimes C_{(1)}$ respectively, see diagram  (\ref{n*}).  
We emphasize that this action is always defined for any quadratic
pair module $C$ and it is natural in the following sense, for any
morphism $f\colon C\r D$ in $\C{qpm}$ and any $n\in\Z$, the equality
$$fn^*=n^*f$$ holds. This property does not hold if $f$ is a weak morphism. The
existence of this action should be compared to the fact that abelian
groups are $\Z$-modules.
\end{defn}


The category $\C{squad}$ of stable quadratic modules is described in \cite{ch4c} IV.C and \cite{2hg1}. Quadratic
modules in general are discussed in \cite{ch4c} and \cite{csch}, they are special $2$-crossed modules in the
sense of \cite{2cm}.
There is a faithful forgetful functor from quadratic pair modules
and weak morphisms to stable quadratic modules
\begin{equation}\label{nse}
\C{wqpm}\To \C{squad}
\end{equation}
sending $C$ as in Definition \ref{qpm} to the stable quadratic
module
\begin{equation}\label{nida} 
\otimes^2(C_{0})_\abb\st{P(-|-)_H}\To C_{1}\st{\partial}\To C_{0}.
\end{equation}
In this paper $G_\abb$ is the abelianization of a group $G$ and $G_\ni$ is its projection of $G$ to the variety of groups of
nilpotency class $2$. 

A track category is a groupoid-enriched category, which is also a $2$-category where all
$2$-morphisms (also termed tracks) are vertically invertible. The category $\C{Top}^*$ of pointed spaces is known to be a track
category with tracks given by homotopy classes of homotopies. The vertical composition in track categories is
denoted by $\vc$, and the vertical inverse of a track $\alpha$ is $\alpha^\vi$. 

The forgetful
functor (\ref{nse}) can be used to pull-back to $\C{wqpm}$ the track
category structure on $\C{squad}$ introduced in \cite{2hg1} 6. The
track structure on $\C{squad}$ was already a pull-back along the
 forgetful functor
\begin{equation}\label{nide}
\C{squad}\To\C{cross}
\end{equation}
from stable quadratic modules to crossed modules considered also in
\cite{2hg1} 6. 

\begin{defn}\label{cros}
We recall that a \emph{crossed module} $\partial\colon M\r N$ is a group homomorphism such that $N$ acts on 
the right of $M$ (the action will be denoted exponentially)
and the homomorphism
$\partial$ satisfies the following two properties $(m,m'\in M, n\in N)$:
\begin{enumerate}
\item $\partial(m^n)=-n+\partial(m)+n$,
\item $m^{\partial(m')}=-m'+m+m'$.
\end{enumerate}
\end{defn}

The crossed module associated via (\ref{nse}) and (\ref{nide}) to a
quadratic pair module $C$ is given by the homomorphism
$$\partial\colon C_{1}\To C_{0},$$ where $C_{0}$ acts on the right
of $C_{1}$ by the formula, $x\in C_{1}$, $y\in C_{0}$,
\begin{equation}\label{nse2}
x^y=x+P(\partial(x)|y)_H.
\end{equation}


\begin{defn}\label{tracks}
A \emph{track} $\alpha\colon f\rr g$ between two morphisms
$f,g\colon C\r D$ in $\C{wqpm}$ is a function
$$\alpha\colon C_{0}\To D_{1}$$
satisfying the equations, $x,y\in C_{0}$, $z\in C_{1}$,
\begin{enumerate}
\item $\alpha(x+y)=\alpha(x)^{f_{0}(y)}+\alpha(y)$,
\item $g_{0}(x)=f_{0}(x)+\partial\alpha(x)$,
\item $g_{1}(z)=f_{1}(z)+\alpha\partial(z)$.
\end{enumerate}
Tracks in $\C{qpm}$ are tracks in $\C{wqpm}$ between morphisms in
the subcategory $\C{qpm}\subset\C{wqpm}$.
\end{defn}

\begin{prop}
The categories $\C{wqpm}$ and $\C{qpm}$ are track categories with
the tracks in Definition \ref{tracks}.
\end{prop}

This proposition is a direct consequence of \cite{2hg1} 6.4.
Vertical and horizontal compositions are defined in the proof of
\cite{2hg1} 6.4.

The following result shows that the weak action of $(\Z,\cdot)$ defined above is also natural with respect to tracks
in $\C{qpm}$.

\begin{prop}
Let $f, g\colon C\r D$ be morphisms in $\C{qpm}$ and let
$\alpha\colon g\rr f$ be a track as in Definition \ref{tracks}. Then
the following diagram commutes
$$\xymatrix{C_0\ar[r]^\alpha\ar[d]_{n^*}&D_1\ar[d]^{n^*}\\C_0\ar[r]^\alpha&D_1}$$
\end{prop}


Given a pointed set $E$ with base point $*\in E$ we denote by $\grupo{E}_\ni$ and $\Z[E]$ to the free group of
nilpotency class $2$ and to the free abelian group generated by $E$ with $*=0$ respectively.

\begin{defn}
A quadratic pair module $C$ is said to be \emph{$0$-free} if
$C_{0}=\grupo{E}_\ni$, $C_{ee}=\otimes^2\Z[E]$
and $H$ is determined by the equalities
$H(e)=0$ for any $e\in E$ and $(s|t)_H=t\otimes s$ for any
$s,t\in \grupo{E}_\ni$. 
\end{defn}

The next lemma shows that $0$-free stable quadratic modules are in
the image of the forgetful functor in (\ref{nse}).

\begin{lem}\label{qpmes}
Any $0$-free stable quadratic module
$$\otimes^2\Z[E]\st{\omega}\To M\st{\partial}\To\grupo{E}_\ni$$
gives rise to a $0$-free quadratic pair module
$$\xymatrix{&\otimes^2\Z[E]\ar[ld]_{P}&\\M\ar[rr]_\partial&&\grupo{E}_\ni\ar[lu]_H}$$
with $P(a\otimes b)=\omega(b\otimes a)$.
\end{lem}


Later we will need the following technical lemma which measures the lack of compatibility of certain tracks in $\C{wqpm}$ with the action of $(\Z,\cdot)$.

\begin{lem}\label{tec}
Let $C$ be a $0$-free quadratic pair module with $C_0=\grupo{E}_\ni$, let  $f\colon C_0\r C_0$ be an endomorphism induced by a pointed map $E\r E$, and let $\alpha\colon C_0\r C_1$ be a map satisfying
\begin{eqnarray*}
\alpha(x+y)&=&\alpha(x)^{f(y)}+\alpha(y),\\
m^*x&=&f(x)+\partial\alpha(x),
\end{eqnarray*}
for some $m\in\Z$ and any $x,y\in C_0$.
Then the following formula holds for any $n\in\Z$ and $x\in C_0$.
$$\alpha(n^*x)=n^*\alpha(x)+\binom{m}{2}\binom{n}{2}P(x|x)_H.$$
\end{lem}

\begin{proof}
We first check that the lemma holds for $x+y$ provided it holds for $x,y\in C_0$. 
\begin{eqnarray*}
\alpha n^*(x+y)&=&\alpha (n^*x+n^*y)\\
&=&\alpha(n^*x)^{f (n^*y)}+\alpha(n^*y)\\
&=&n^*\alpha(x)+n^*\alpha(y)+\binom{m}{2}\binom{n}{2}P(x|x)_H\\
&&\binom{m}{2}\binom{n}{2}P(y|y)_H+P(-f (n^*x)+n^*m^*x|f (n^*y))_H\\
&=&n^*(\alpha(x)+\alpha(y))+n^2P(-f (x)+m^*x|f y)_H
\\&&+
\binom{m}{2}\binom{n}{2}P(x+y|x+y)_H\\
&=&n^*(\alpha(x)^{f (y)}+\alpha(y))+
\binom{m}{2}\binom{n}{2}P(x+y|x+y)_H\\
&=&n^*\alpha(x+y)+
\binom{m}{2}\binom{n}{2}P(x+y|x+y)_H.
\end{eqnarray*}
Here we use that $f$ is compatible with the action of $(\Z,\cdot)$ and that $P(x|x)_H$ is linear in $x$.

Now since $C_0=\grupo{E}_\ni$ we only need to check that the proposition holds for $e\in E$. But $H(e)=0$, so we have $n^*e=n\cdot e$. The equality
\begin{eqnarray*}
\alpha(n\cdot e)  
&=&n\cdot\alpha(e)+\binom{n}{2}P(f (e)|f (e))_H+m\binom{n}{2}P(e|f (e))_H
\end{eqnarray*}
follows easily by induction in $n$ from the first equation of the statement and the laws of a quadratic pair module. On the other hand 
\begin{eqnarray*}
n^*\alpha(e)&=&n\cdot\alpha(e)+\binom{n}{2}PH\partial(-f (e)+m\cdot e).
\end{eqnarray*}
One can also check by induction that
\begin{eqnarray*}
PH(-f (e)+m\cdot e)&=&P(f (e)|f (e))_H+mP(e|f (e))_H-\binom{m}{2}P(e|e)_H.
\end{eqnarray*}
Now the proof is finished.
\end{proof}

Lemma \ref{tec} holds under the more general condition that $C_0$ is generated by elements $x\in C_0$ with $H(x)=0$
and $Hf(x)=0$.

\section{Homotopy groups and secondary homotopy groups}\label{el1}

Let $\C{Top}^*$ be the category of (compactly generated)
pointed spaces. Using classical homotopy groups $\pi_nX$ we obtain for $n\geq 0$ the functor
$$\Pi_n\colon\C{Top}^*\To\C{Ab}$$
with
\begin{equation}\label{Clas}
\Pi_nX=\left\{
\begin{array}{ll}
\pi_nX, & n\geq2,\\
(\pi_1X)_\abb, & n=1,\\
\Z[\pi_0X], & n=0,
\end{array}\right.
\end{equation}
termed \emph{additive homotopy group}. Here $G_\abb$ is the abelianization of a group $G$. 

One readily checks that the smash product
$$f\wedge g\colon S^n\wedge S^m\To X\wedge Y$$
of maps $\set{f\colon S^n\r X}\in\pi_nX$ and $\set{g\colon S^m\r Y}\in\pi_mY$
induces a well-defined homomorphism
\begin{equation}\label{ClasS}
\wedge\colon\Pi_nX\otimes\Pi_mY\To\Pi_{n+m}(X\wedge Y).
\end{equation}
This homomorphism is symmetric in the sense that the interchange map $\tau_{X,Y}\colon X\wedge Y\r Y\wedge
X$ yields the equation in $\Pi_{n+m}(Y\wedge X)$
\begin{equation}\label{inter}
(\tau_{X,Y})_*(f\wedge g)=(-1)^{nm}g\wedge f.
\end{equation}
Here the sign $(-1)^{nm}$ is given by the interchange map
\begin{equation}\label{inter2}
\tau_{n,m}=\tau_{S^n,S^m}\colon S^{n+m}\To S^{m+n}
\end{equation}
which has degree $(-1)^{nm}$. Here $\tau_{n,m}$ also designates the corresponding element of the symmetric group $\sym{n+m}$
which acts from the left on $S^{n+m}$, see Section \ref{the} below.

We want to generalize the smash product operator (\ref{ClasS}) for 
additive secondary
homotopy groups. 

\begin{defn}
Let $n\geq 2$. For a pointed space $X$ we define the \emph{additive secondary homotopy group} $\Pi_{n,*}X$ which
is the $0$-free quadratic pair module given by the diagram 
$$\Pi_{n,*}X=\left(\begin{array}{c}\xymatrix@C=0pt{&**[r]\Pi_{n,ee}X=\otimes^2\Z[\L^nX]\ar[ld]_{P}&\\\Pi_{n,1}X\ar[rr]_\partial&&**[r]\Pi_{n,0}X=\grupo{\L^nX}_\ni\ar[lu]_H}\end{array}\right)$$
We obtain the group $\Pi_{n,1}X$ and the homomorphisms $P$ and
$\partial$ as follows. The group $\Pi_{n,1}X$  is given by the set of equivalence classes $[f,F]$ 
represented
by a map $f\colon S^1\r \vee_{\L^n X}S^1$ and a track
$$\xymatrix{S^n\ar[r]_{\S^{n-1}f}^<(.98){\;\;\;\;\;}="a"\ar@/^25pt/[rr]^0_{}="b"&S^n_X\ar[r]_{ev}&X.\ar@{=>}"a";"b"_F}$$
Here the pointed space $$S^n_X=\vee_{\L^n X}S^n=\S^n\L^n X$$ is the $n$-fold suspension of the $n$-fold loop space
$\L^nX$, where $\L^nX$ is regarded as a pointed set with the discrete topology. Hence $S^n_X$ is
the coproduct of $n$-spheres indexed by the set of non-trivial maps $S^n\r X$, and $ev\colon S^n_X\r X$
is the obvious evaluation map. Moreover, for the sake of simplicity given a map
$f\colon S^1\r \vee_{\L^n X}S^1$ we will denote $f_{ev}=ev(\S^{n-1}f)$, so that $F$ in the previous diagram is a track $F\colon f_{ev}\rr 0$.
The equivalence relation $[f,F]=[g,G]$ holds provided there is a 
track $N\colon\Sigma^{n-1}f\rr\S^{n-1}g$ with $\hopf(N)=0$ if $n\geq 3$ or $\bar{\sigma}\hopf(N)=0$ if
$n=2$, see (\ref{elhopf}) and (\ref{sibar}) below, such that 
the composite track in the following diagram is the trivial track.
$$\xymatrix@C=50pt{S^n\ar@/^40pt/[rr]^0_{\;}="a"\ar@/_40pt/[rr]_0^{\;}="f"\ar@/^15pt/[r]|{\S^{n-1}f}^<(.935){\;}="b"_{\;}="c"\ar@/_15pt/[r]|{\S^{n-1}g}^{\;}="d"_<(.93){\;}="e"
&S^n_X\ar[r]^{ev}&X\ar@{=>}"a";"b"^{F^\vi}\ar@{=>}"c";"d"^{N_{f,g}}\ar@{=>}"e";"f"^G}$$
That is $F=G\vc(ev\,N_{f,g})$. The map
$\partial$ is defined by the  formula
$$\partial[f,F]=(\pi_1 f)_\ni(1),$$
where $1\in\pi_1S^1=\Z$.

The Hopf invariant of a track $N\colon\Sigma^{n-1}f\rr\S^{n-1}g$ as above is defined in \cite{2hg1} 3.3 by the
homomorphism 
\begin{equation}\label{elhopf}
H_2(IS^1,S^1\vee S^1)\st{ad(N)_*}\To H_2(\L^{n-1}S^n_X,\vee_{\L^ nX}S^1)\cong\left\{
\begin{array}{ll}
\hat{\otimes}^2\Z[\L^nX],& n\geq 3,\\&\\
{\otimes}^2\Z[\L^nX],& n=2.
\end{array}
\right.
\end{equation}
which carries the generator $1\in\Z\cong H_2(IS^1,S^1\vee S^1)$ to $\hopf(N)$. Here $ad(N)_*$ is the homomorphism
induced in homology by the adjoint of the homotopy $$N\colon \S^{n-1} IS^1\cong IS^n\r S^n_X.$$ The
reduced tensor square is given by
$$\hat{\otimes}^2A=\frac{A\otimes A}{a\otimes b+b\otimes a\sim 0},$$ 
and 
\begin{equation}\label{sibar}
\bar{\sigma}\colon\otimes^2A\twoheadrightarrow\hat{\otimes}^2A
\end{equation}
is the natural projection. The isomorphism in (\ref{elhopf}) is induced by the Pontrjaging product. 
We refer the reader to \cite{2hg1} 3 for a
complete definition of the Hopf invariant for tracks and for the elementary properties which will be used in this paper.
For the sake of simplicity we define the \emph{reduced Hopf invariant} as $\overline{\hopf}=\hopf$ if $n\geq 3$
and
$\overline{\hopf}=\bar{\sigma}\hopf$ if $n=2$. A \emph{nil-track} in this paper will be a track in $\C{Top}^*$ with trivial
reduced Hopf invariant. In particular the equivalence relation defining elements in $\Pi_{n,1}X$ is determined by
nil-tracks.

This completes the definition
of $\Pi_{n,1}X$, $n\geq 2$, as a set. The
group structure of $\Pi_{n,1}X$ is induced by the comultiplication $\mu\colon S^1\r S^1\vee S^1$, compare
\cite{2hg1} 4.4.

We now define the homomorphism $P$ for additive secondary homotopy groups $\Pi_{n,*}X$ with $n\geq 2$. Consider the diagram
$$\xymatrix{S^n\ar[rr]_{\S^{n-1}\beta}^{\;}="a"\ar@/^30pt/[rr]^0_{\;}="b"&&S^n\vee S^n\ar@{=>}"a";"b"_{B}}$$
where $\beta\colon S^1\r S^1\vee S^1$ is given such that $(\pi_1\beta)_\ni(1)=-a-b+a+b\in\grupo{a,b}_\ni$ is the commutator. The track $B$ 
is any track with
$\overline{\hopf}(B)=-\sigma(a\otimes b) \in\hat{\otimes}^2\Z[a,b]$. 
Given $x\otimes y\in\otimes^2\Z[\L^nX]$ let
$\tilde{x},\tilde{y}\colon S^1\r\vee_{\L^nX}S^1$ be maps with
$(\pi_1\tilde{x})_\abb(1)=x$ and $(\pi_1\tilde{y})_\abb(1)=y$. Then
the diagram
\begin{equation}\label{omedia}
\xymatrix{S^n\ar[rr]_{\S^{n-1}\beta}^{\;}="a"\ar@/^30pt/[rr]^0_{\;}="b"&&S^n\vee
S^n\ar@{=>}"a";"b"_{B}\ar[rr]_{\S^{n-1}(\tilde{y},\tilde{x})}&&S^n_X\ar[r]_{ev}&X}
\end{equation}
represents an element
\begin{equation*}
P(x\otimes y)=[(\tilde{y},\tilde{x})\beta,ev\,(\S^{n-1}(\tilde{y},\tilde{x}))B]\in\Pi_{n,1}X.
\end{equation*}
This completes the definition of the quadratic pair module $\Pi_{n,*}X$ for $n\geq 2$. For $n=0,1$ we define the
additive secondary homotopy groups $\Pi_{n,*}X$ by the following remark. In this way we get for $n\geq 0$ a
functor
$$\Pi_{n,*}\colon\C{Top}^*\To\C{qpm}$$
which is actually a track functor.
\end{defn}

\begin{rem}
Considering maps $f\colon S^n\r X$ together with tracks of
such maps to the trivial map, we introduced in \cite{2hg1} the secondary
homotopy group $\pi_{n,*}X$, which is a groupoid for $n=0$, a crossed module for $n=1$, a reduced quadratic
module for $n=2$, and a stable quadratic module for $n\geq 3$. Let $\C{squad}$ be the category of stable quadratic
modules. 

Then using the adjoint functors $\ad_n$ of the forgetful functors $\phi_n$ as discussed in \cite{2hg1} 6 we get the 
\emph{additive secondary homotopy group} track functor
\begin{equation*}
\Pi_{n,*}\colon\C{Top}^*\To\C{squad}
\end{equation*}
given by
$$\Pi_{n,*}X=\left\{\begin{array}{ll}
\pi_{n,*}X,&\text{for }n\geq3,\\
\ad_3\pi_{2,*}X,&\text{for }n=2,\\
\ad_3\ad_2\pi_{1,*}X,&\text{for }n=1,\\
\ad_3\ad_2\ad_1\pi_{0,*}X,&\text{for }n=0.\\
\end{array}\right.$$
This is the secondary analogue of (\ref{Clas}).

Here the category $\C{squad}$ of stable quadratic modules is not appropriate to study the smash product of
secondary homotopy groups since we do not have a symmetric monoidal structure in $\C{squad}$. Therefore 
we introduced above the category
$\C{qpm}$ of quadratic pair modules and we observe that $\Pi_{n,*}X$ in $\C{squad}$ yields a functor
to the category $\C{qpm}$ as follows.
A map $f\colon X\r Y$ in $\C{Top}^*$ induces a homomorphism
$\Pi_{n,0}f\colon\Pi_{n,0}X\r\Pi_{n,0}Y$ between free nil-groups
which carries generators in $\Pi_{n,0}X$ to generators in
$\Pi_{n,0}Y$ and therefore $\Pi_{n,*}f$ is compatible with $H$. This shows that Lemma \ref{qpmes} gives rise to a canonical lift
\begin{equation*}
\xymatrix{&\C{qpm}\ar[d]\\\C{Top}^*\ar[ru]^{\Pi_{n,*}}\ar[r]_{\Pi_{n,*}}&\C{squad}}
\end{equation*}
Here the vertical arrow, which is the forgetful track functor given by (\ref{nida}), is faithful but
not full at the level of morphisms. 

The definition of $\Pi_{2,*}X$ given above coincides with the lifting of $\ad_3\pi_{2,*}X$ to $\C{qpm}$ by the
claim (*) in the proof of \cite{2hg1} 4.9.
\end{rem}

In this paper we are concerned with the properties of the
track functor $\Pi_{n,*}$, mapping to the category $\C{qpm}$.
The category $\C{qpm}$ is, in fact, a symmetric monoidal category, defined by a tensor product $\odot$ in
$\C{qpm}$, see \cite{qaI}, and the smash product yields the operator
\begin{equation}
\wedge\colon\Pi_{n,*}X\odot\Pi_{m,*}Y\To\Pi_{n+m,*}(X\wedge Y)
\end{equation}
constructed in \cite{2hg3}. 
Equation (\ref{inter}) has now a secondary analogue
given by the right action of the symmetric group $\sym{n+m}$ on the object
$\Pi_{n+m,*}(X\wedge Y)$ in $\C{qpm}$. More precisely the following diagram commutes in $\C{qpm}$.
$$\xymatrix{\Pi_{n,*}X\odot\Pi_{m,*}Y\ar[dd]^\cong_{\tau_\odot}\ar[r]^{\wedge}&\Pi_{n+m,*}(X\wedge
Y)\ar[d]_\cong^{(\tau_{X,Y})_*}\\&\Pi_{n+m,*}(Y\wedge X)\ar@{<-}[d]_\cong^{\tau_{n,m}^*}\\
\Pi_{m,*}Y\odot\Pi_{n,*}X\ar[r]^{\wedge}&\Pi_{m+n,*}(Y\wedge X)}$$
Here $\tau_\odot$ on the left hand side is given by the symmetry of the tensor product $\odot$ in $\C{qpm}$ and
$\tau_{n,m}^*$ is defined by the action of $\sym{n+m}$.
For this reason we define and 
study in this paper the properties of the symmetric
group action on secondary homotopy groups.




\section{Actions of monoid-groupoids in track categories}

In this paper we deal with actions on additive secondary homotopy groups.
Additive secondary homotopy groups are objects in a track category. In
ordinary categories a monoid action is given by a
monoid-morphism mapping to an endomorphism monoid in the category.
In track categories endomorphism objects are monoids in the monoidal category
of groupoids, where the
monoidal structure is given by the (cartesian) product. Therefore one can define
accordingly actions of such monoids. We make explicit this structure in the following
definition.

\begin{defn}\label{mg}
Let $\mathbb{I}$ be the
category with only one object $*$ and one morphism $1\colon *\r *$.
A \emph{monoid-groupoid} $\C{G}$ is a groupoid together with a
multiplication functor $\cdot\colon\C{G}\times\C{G}\r \C{G}$ and a
unit functor $u\colon\mathbb{I}\r \C{G}$, satisfying the laws
of a monoid in the symmetric monoidal category of groupoids. We usually identify $*=u(*)$. The opposite
$\C{G}^{op}$ of a monoid-groupoid is the  underlying groupoid
$\C{G}$ with its unit functor and multiplication functor given
by
$$\C{G}\times\C{G}\st{T}\To\C{G}\times\C{G}\st{\cdot}\To\C{G}.$$
Here $T$ is the interchange of factors in the product. 
A monoid-groupoid \emph{morphism} $f\colon\C{G}\r\C{H}$ is a functor
preserving the multiplication and the unit.
\end{defn}

Monoid-groupoids are also termed strict monoidal groupoids. The
weaker versions of this concept will not be considered in this
paper, therefore we abbreviate the terminology.

The canonical example of a monoid-groupoid is obtained by the endomorphisms of
an object $X$ in a track category $\C{C}$, denoted
$$\enda_\C{C}(X).$$
The multiplication is given by composition in $\C{C}$, and the unit
is given by the identity morphism $1_X\colon X\r X$. In fact a
monoid-groupoid as defined above is exactly the same thing as a
track category with only one object, the opposite monoid-groupoid
coincides with the the opposite of the corresponding track category
and monoid-groupoid morphisms correspond to $2$-functors.

\begin{defn}\label{mga}
Let $X$ be an object in a track category $\C{C}$ and let $\C{G}$ be
a monoid-groupoid. A \emph{right action} of $\C{G}$ on $X$ is a
monoid-groupoid morphism $\C{G}^{op}\r\enda_\C{C}(X)$.
\end{defn}

An important example of monoid-groupoids arises from crossed
modules. The monoid-groupoid $M(\partial)$ associated to a
crossed module $\partial\colon T\r G$ has object set $G$ and
morphism set the semidirect product $G\ltimes T$.
Here we write the groups $T$ and $G$ with a multiplicative group
law. An element $(g,t)\in G\ltimes T$ is a morphism $(g,t)\colon
g\cdot\partial(t)\r g$ in $M(\partial)$. The composition law $\circ$ is given by the formula
$(g,t)\circ(g\cdot\partial(t),t')=(g,t\cdot t')$. Multiplication in the groups $G$
and $G\ltimes T$ defines the multiplication of $M(\partial)$ and
the unit is given by the unit elements in $G$ and $G\ltimes T$.
Indeed this correspondence determines an equivalence between crossed
modules and group objects in the category of groupoids. This example
can be used to define crossed module actions.

\begin{defn}\label{cma}
Let $X$ be an object in a track category $\C{C}$ and let
$\partial\colon T\r G$ be a crossed module. A \emph{right action} of
$\partial\colon T\r G$ on $X$ is a monoid-groupoid morphism
$M(\partial)^{op}\r\enda_\C{C}(X)$. 
\end{defn}

We are interested in right actions of crossed modules in the track
category $\C{wqpm}$. We explicitly describe such actions as follows.

\begin{defn}\label{cma2}
A \emph{right action of a crossed module} $\partial\colon T\r G$ on
a quadratic pair module $C$ in the category $\C{wqpm}$ consists of a
group action of $G$ on the right of $C$ given by morphisms in
$\C{wqpm}$,
$$g^*\colon C\To C,\;\;g\in G,$$
together with a bracket
$$\grupoo{-,-}\colon C_0\times T\To C_1$$
satisfying the following properties, $x,y\in C_0$, $z\in C_1$,
$s,t\in T$, $g\in G$,
\begin{enumerate}
\item $\grupoo{x+y,t}=\grupoo{x,t}^{\partial(t)^*y}+\grupoo{y,t}$,
\item $x=\partial(t)^*x+\partial\grupoo{x,t}$,
\item $z=\partial(t)^*z+\grupoo{\partial(z),t}$,
\item $\grupoo{x,s\cdot t}=\grupoo{\partial(s)^*x,t}+\grupoo{x,s}=\partial(t)^*\grupoo{x,s}+\grupoo{x,t}$,
\item $\grupoo{x,t^g}=g^*\grupoo{(g^{-1})^*x,t}$.
\end{enumerate}
We point out that the second equality in (4) follows from (1)--(3).
Indeed these are the two possible definitions of the horizontal
composition $\grupoo{-,s}\grupoo{-,t}\colon\partial(st)^*\rr 1$ of the tracks $\grupoo{-,t}\colon\partial(t)^*\rr 1$ and
$\grupoo{-,s}\colon\partial(s)^*\rr 1$ in the track category
$\C{wqpm}$.
\end{defn}

The notion of action defined above corresponds to an action in
Norrie's sense (\cite{norrie}) of a crossed module on the underlying
crossed module of a quadratic pair
module, however Norrie considers left actions. 

The very special kind of action introduced in the following
definition will be of importance to describe the symmetric action on additive 
secondary homotopy groups in Section \ref{the}.

\begin{defn}\label{signg}
Let $\set{\pm1}$ be the multiplicative group of order $2$. A
\emph{sign group} $G_\vc$ is a diagram of group homomorphisms
$$\set{\pm1}\st{\imath}\hookrightarrow G_\vc\st{\partial}\twoheadrightarrow G\st{\varepsilon}\To \set{\pm1}$$
where the first two morphisms form an extension. Here all groups
have a multiplicative group law and the composite
$\varepsilon\partial$ is also denoted by $\varepsilon\colon G_\vc\r
\set{\pm1}$. Moreover, we define the element $\omega=\imath(-1)\in G_\vc$.

A sign group $G_\vc$ \emph{acts on the right of a quadratic pair module}
$C$ if $G$ acts on the right of $C$ by morphisms $$g^*\colon C\To
C,\;\; g\in G,\;\;\text{ in }\C{qpm},$$ and there is a bracket
$$\grupo{-,-}\colon C_0\times G_\vc\To C_1$$
satisfying the following properties, $x,y\in C_0$, $z\in C_1$,
$s,t\in G_\vc$, were $\varepsilon(t)^*$ is given by the action of
$(\Z,\cdot)$ in Definition \ref{sqmact},
\begin{enumerate}
\item $\grupo{x+y,t}=\grupo{x,t}^{\partial(t)^*y}+\grupo{y,t}$,
\item $\varepsilon(t)^*(x)=\partial(t)^*(x)+\partial\grupo{x,t}$,
\item $\varepsilon(t)^*(z)=\partial(t)^*(z)+\grupo{\partial(z),t}$,
\item $\grupo{x,s\cdot
t}=\grupo{\partial(s)^*(x),t}+\grupo{\varepsilon(t)^*x,s}$,
\item the \emph{$\omega$-formula}:
$$\grupo{x,\omega}=P(x|x)_H.$$
\end{enumerate}
Notice that the $\omega$-formula corresponds to the $k$-invariant,
see \cite{2hg1} 8.
\end{defn}

\begin{rem}
A sign group $G_\vc$ gives rise to a crossed module
$$\partial_\vc=(\varepsilon,\partial)\colon G_\vc\To\set{\pm1}\times G,$$
where $\set{\pm1}\times G$ acts on
$G_\vc$ by the formula
$$g^{(x,h)}=\bar{h}^{-1}g\bar{h}\imath\left(\varepsilon(g)^{\binom{x}{2}}\right).$$
Here $g\in G_\vc$, $x\in\set{\pm1}$, $h\in G$ and $\bar{h}\in G_\vc$ is any element
with $\partial(\bar{h})=h$. This action is well defined since
$G_\vc$ is a central extension of $G$ by $\set{\pm1}$.
\end{rem}

\begin{lem}
The sign group action in Definition \ref{signg} corresponds to an
action of the crossed module $\partial_\vc$ on $C$ in the sense of
Definition \ref{cma2} such that $\set{\pm1}$ acts on $C$ by the
action of $(\Z,\cdot)$ in Definition \ref{sqmact}, $G$ acts by
morphisms in $\C{qpm}$, and the $\omega$-formula holds. The
correspondence is given by the formula
$$\grupoo{x,t}=\grupo{\varepsilon(t)^*x,t},\;\;x\in C_0,\;\;t\in G_\vc.$$
\end{lem}

The proof of this lemma is straightforward. We just want to point out that Definition \ref{cma2} (5) follows in
this case from Definition \ref{signg} (4), (5), and Lemma \ref{tec}. 

\begin{rem}\label{pursu}
A sign group $G_\vc$ is \emph{trivial} if $G$ is a trivial group. Notice that a trivial sign group acts on
any quadratic pair module in a unique way.
\end{rem}

\section{The action of $\endo_*(S^n)$ on $\Pi_{n,*}X$}\label{todas}

Let $S^n$ be the $n$-sphere and let $\endo_*(S^n)=\L^n
S^n$ be the topological monoid of maps $S^n\r S^n$ in $\C{Top}^*$. Then the fundamental groupoid of
$\endo_*(S^n)$, denoted by $\pi_{0,*}\endo_*(S^n)$, is a monoid-groupoid in the sense of Definition \ref{mg}.
It is well known that the monoid of path components of
$\endo_*(S^n)$ coincides with the multiplicative monoid
$(\Z,\cdot)$. 

We now consider the right action of $\pi_{0,*}\endo_*(S^n)$ on
$\Pi_{n,*}X$ for $n\geq 2$. That is, we define for each pointed map
$f\colon S^n\r S^n$ an induced map in $\C{qpm}$
$$f^*\colon\Pi_{n,*}X\To\Pi_{n,*}Y$$
and we define for each track $H\colon f\rr g$ with $f,g\colon S^n\r
S^n$ a track in $\C{qpm}$
$$H^*\colon f^*\rr g^*.$$
This yields a right action of the fundamental groupoid
$\pi_{0,*}\endo_*(S^n)$ on the secondary homotopy group $\Pi_{n,*}X$
in the track category $\C{qpm}$ of quadratic pair modules in the
sense of Definition \ref{mga}. 

\begin{thm}\label{laa}
Let $X$ be a pointed space. For any $n\geq 2$ there is a natural
action of the monoid-groupoid $\pi_{0,*}\endo_*(S^n)$ on the
quadratic pair module $\Pi_{n,*}X$.
\end{thm}

The rest of this section is devoted to the proof of this theorem,
which is carried out in several steps. 

The discrete monoid $\pi_{0,0}\endo_*(S^n)$, which is the underlying 
set of the topological monoid $\endo_*(S^n)$, acts on the right of
the pointed set $\L^n X$ of pointed maps $S^n\r X$ by
precomposition, i. e. given $f\colon S^n\r S^n$ the induced
endomorphism is
$$f^*\colon\L^nX\To\L^nX,\;\; f^*(g)=gf.$$
This induces a right action of $\pi_{0,0}\endo_*(S^n)$
on the free group $\pi_{n,*}X=\grupo{\L^nX}_\ni$ of nilpotency class $2$ which will be
denoted in the same way.

In order to extend this action to $\Pi_{n,1}X$ we consider the
submonoid
\begin{equation}\label{nse4}
\tilde{\pi}_{0,1}\endo_*(S^n)\subset\pi_{0,1}\endo_*(S^n)
\end{equation}
of the monoid $\pi_{0,1}\endo_*(S^n)$ of morphisms in $\pi_{0,*}\endo_*(S^n)$
given by tracks between self-maps of $S^n$ of the form
\begin{equation}\label{nseya}
\gamma\colon f\rr \S^{n-1}(\cdot)^{\deg f}=(\cdot)^{\deg f}_n.
\end{equation}
Here $\deg f\in\Z$ denotes the degree of $f\colon S^n\r S^n$ and for $k\in\Z$
$$(\cdot)^k\colon S^1\To S^1\colon z\mapsto z^k$$
is given by the (multiplicative) topological abelian group structure
of $S^1$. 

We need a bracket operation
\begin{equation}\label{bra}
\grupo{-,-}\colon \Pi_{n,0}X\times
\tilde{\pi}_{0,1}\endo_*(S^n)\To\Pi_{n,1}X,
\end{equation}
defined as follows. Let $x\in\pi_{n,0}X=\grupo{\L^nX}_\ni$ and $\gamma\colon
f\rr(\cdot)^{\deg f}_n$ in $\tilde{\pi}_{0,1}\endo_*(S^n)$. We
choose maps $\tilde{x}\colon S^1\r\vee_{\L^nX}S^1$,
$\varepsilon\colon S^1\r S^1\vee S^1$ with
$(\pi_1\tilde{x})_\ni(1)=x$ and $(\pi_1\varepsilon)_\ni=-a+b\in\grupo{a,b}_\ni$. Then
$\grupo{x,\gamma}\in\Pi_{n,1}X$ is the element represented by the
map
\begin{equation*}
\xymatrix{S^1\ar[r]^<(.25)\varepsilon& S^1\vee
S^1\ar[r]^<(.15){\tilde{x}\vee\tilde{x}}&(\vee_{\L^nX}S^1)\vee(\vee_{\L^nX}S^1)
\ar[rrr]^<(.4){(\S f^*,\vee_{\L^nX}(\cdot)^{\deg f})}
&&&\vee_{\L^nX}S^1}
\end{equation*}
and the track
\begin{equation}\label{detra}
\xymatrix@C=35pt@R=60pt{&S^n\ar[r]^{\tilde{x}}&S^n_X\ar[rrr]^{\vee_{\L^nX}(\cdot)_n^{\deg
f}}_<(.4){\;}="f" &&&S^n_X\ar@/^20pt/[rd]^{ev}&
\\S^n\ar@/^20pt/[ru]^0_{\;}="b"\ar[r]_{\S^{n-1}\varepsilon}&
S^n\vee
S^n\ar[u]_{(1,1)}^<(.1){\;}="a"\ar[r]_{\S^{n-1}(\tilde{x}\vee\tilde{x})}&
S^n_X\vee S^n_X\ar[u]^{(1,1)}
\ar[rrru]_<(.3){\;\;\;\;\;\;\;\;\;(\vee_{\L^nX} f,\vee_{\L^nX}(\cdot)_n^{\deg
f})}^<(.5){\;}="e"
\ar[rrr]_{(\S^n f^*,\vee_{\L^nX}(\cdot)_n^{\deg f})}
&&&S^n_X\ar[r]_{ev}&X\ar@{=>}"a";"b"^N\ar@{=>}"e";"f"^{(\vee_{\L^nX}\gamma,0^\vc)}}
\end{equation}
Here $N$ is a nil-track.

The main properties of the bracket operation in (\ref{bra}) are
listed in the following proposition.

\begin{prop}\label{nores}
The bracket $\grupo{-,-}$ in (\ref{bra}) satisfies the following
formulas for any $x,y\in\Pi_{n,0}X$ and $\gamma\colon
f\rr(\cdot)_n^{\deg f},\delta\colon g\rr(\cdot)_n^{\deg g}$ in
$\tilde{\pi}_{0,1}\endo_*(S^n)$.
\begin{enumerate}
\item $\grupo{x+y,\gamma}=\grupo{x,\gamma}^{f^*y}+\grupo{y,\gamma}$,
\item $(\deg f)^*x=f^*x+\partial\grupo{x,\gamma}$,
\item $\grupo{x,\gamma\delta}=\grupo{f^*x,\delta}+\grupo{(\deg
g)^*x,\gamma}$,
\item if $\omega\colon 1_{S^n}\rr 1_{S^n}$ is a track with $0\neq\overline{\hopf}(\omega)\in\hat{\otimes}^2\Z=\Z/2$
then $\grupo{x,\omega}=P(x|x)_H$.
\end{enumerate}
Moreover, this bracket operation is natural in $X$.
\end{prop}

\begin{proof}
With the notation in \cite{2hg1} 7.4 we have
$\grupo{x,\gamma}=r(ev(\vee_{\L^nX}\gamma)(\S^{n-1}\tilde{x}))$ for
the track
\begin{small}
$$ev(\vee_{\L^nX}\gamma)(\S^{n-1}\tilde{x})\colon ev(\S^nf^*)(\S^{n-1}\tilde{x})=
ev(\vee_{\L^nX}f)(\S^{n-1}\tilde{x})\rr ev(\vee_{\L^nX}(\cdot)_n^{\deg
f})(\S^{n-1}\tilde{x}),$$ 
\end{small}
therefore (1) and (2) follow from
\cite{2hg1} 7.6 and 7.5 (2).

It is easy to see that the formula
$$ev(\vee_{\L^nX}\gamma\delta)=(ev(\vee_{\L^nX}\gamma)(\vee_{\L^nX}(\cdot)^{\deg g}_n))\vc(ev(\vee_{\L^nX}\delta)(\S^nf^*))$$
holds, therefore (3) follows from \cite{2hg1} 7.5 (3).


If we evaluate $\grupo{x,-}$ at $\omega$ then the composite track obtained from (\ref{detra}) 
by going from the lower left $S^n$ to the upper right $S^n_X$ has the same reduced Hopf invariant as the 
track from $S^n$ to $S^n_X$ in (\ref{omedia}). Indeed the formula for both reduced Hopf invariants is (c) 
in the proof of Proposition \ref{oes}. Therefore (4) follows.
\end{proof}

The next result follows from the algebraic properties of the
bracket (\ref{bra}) which are proved in the previous proposition together with Lemma \ref{tec}.

\begin{prop}\label{act}
The monoid $\tilde{\pi}_{0,1}\endo_*(S^n)$ acts on the right of
$\Pi_{n,1}X$ by the following formula, $n\geq 2$: given
$x\in\Pi_{n,1}X$ and $\gamma\colon f\rr(\cdot)_n^{\deg f}$
$$\gamma^*x=(\deg f)^*x-\grupo{\partial(x),\gamma}.$$
This action satisfies $\partial\gamma^*=f^*\partial$,
$\gamma^*P=P(\otimes^2 f^*_\abb)$, and $Hf^*=(\otimes^2f^*_\abb)H$, therefore it defines an
action of $\tilde{\pi}_{0,1}\endo_*(S^n)$ on the right of the
quadratic pair module $\Pi_{n,*}X$ in the category $\C{qpm}$. 
This action is natural in $X$.
\end{prop}

\begin{proof}
The equality $Hf^*=(\otimes^2f^*_\abb)H$ follows from the fact that the endomorphism $f^*$ 
carries generators to generators in $\grupo{\L^nX}_\ni$. 
The equality $\partial\gamma^*=f^*\partial$ follows from Proposition \ref{nores} (2). Let us check $\gamma^*P=P(\otimes^2 f^*_\abb)$. Given $a,b\in \L^nX$
\begin{eqnarray*}
\gamma^*P(a\otimes b)&=&(\deg f)^*P(a\otimes b)-\grupo{-a-b+a+b,\gamma}\\
&=&P(\deg f)^2(a\otimes b)-\grupo{b,\gamma}-\grupo{a,\gamma}+\grupo{b,\gamma}+\grupo{a,\gamma}\\
&&+P(-f^*(a)+(\deg f)^*a|f^*b)_H-P(-f^*(b)+(\deg f)^*b|f^*a)_H\\
&=&P(\deg f)^2(a\otimes b)+P(\partial\grupo{b,\gamma}|\partial\grupo{a,\gamma})_H\\
&&+P(-f^*(a)+(\deg f)^*a|f^*b)_H-P(-f^*(b)+(\deg f)^*b|f^*a)_H\\
&=&-P(\deg f)^2(a| b)_H\\&&-P(-f^*(a)+(\deg f)^*a|-f^*(b)+(\deg f)^*b)_H\\
&&+P(-f^*(a)+(\deg f)^*a|f^*b)_H+P(f^*a|-f^*(b)+(\deg f)^*b)_H\\
&=&-P(f^*(a)|f^*(b))_H\\
&=&P(f^*(b)|f^*(a))_H\\
&=&P(f^*(a)\otimes f^*(a)).
\end{eqnarray*}
Here we use Proposition \ref{nores} (1) and (2) and the fact that $H(a)=0=H(b)$.

Finally given $\delta\colon g\rr(\cdot)_n^{\deg g}$
\begin{eqnarray*}
\gamma^*\delta^*(x)&=&\gamma^*((\deg g)^*x-\grupo{\partial(x),\delta})\\
&=&(\deg f)^*(\deg g)^*x-(\deg f)^*\grupo{\partial(x),\delta}
-\grupo{g^*\partial(x),\gamma}\\
&=&((\deg f)(\deg g))^*x-\grupo{(\deg f)^*\partial(x),\delta}
-\grupo{g^*\partial(x),\gamma}
\\&&+\binom{\deg f}{2}\binom{\deg g}{2}P(\partial(x)|\partial(x))_H\\
&=&(\deg fg)^*x-\grupo{\partial(x),\gamma\delta}\\
&=&(\gamma\delta)^*(x).
\end{eqnarray*}
Here we use Proposition \ref{nores} (1), (2) and (3), Lemma \ref{tec} and the fact that 
$P(\partial(x)|\partial(x))_H=-x-x+x+x=0$.
\end{proof}

\begin{prop}\label{oes}
For $n\geq 2$ the right action of the monoid $\tilde{\pi}_{0,1}\endo_*(S^n)$ on the
group  $\Pi_{n,1}X$ given by Proposition
\ref{act} factors through the boundary homomorphism
$$
q\colon\tilde{\pi}_{0,1}\endo_*(S^n)\twoheadrightarrow
\pi_{0,0}\endo_*(S^n),\;\; q(\gamma\colon f\rr(\cdot)_n^{\deg
f})=f,$$
that is, the homomorphism $\gamma^*=f^*$ only depends on the boundary $q(\gamma)=f$.
\end{prop}

\begin{proof}
Let $\gamma\colon f\rr(\cdot)^{\deg f}_n$ be any element in
$\tilde{\pi}_{0,1}\endo_*(S^n)$ and let $\delta\colon (\cdot)^{\deg
f}_n\rr (\cdot)^{\deg f}_n$ be any track.
We know that all elements in $q^{-1}(f)$ are of the form $\delta\vc\gamma$,
therefore we only have to check that for any $[g,G]\in\Pi_{n,1}X$
$$\gamma^*[g,G]=(\delta\vc\gamma)^*[g,F],$$
or equivalently
$$\grupo{\partial[g,G],\gamma}=\grupo{\partial[g,G],\delta\vc\gamma}.$$


The element $\grupo{\partial[g,G],\delta\vc\gamma}$ is represented
by the following diagram
\begin{equation*}\tag{a}
\xymatrix@R=100pt{&S^n\ar[r]^{\S^{n-1}g}&S^n_X\ar[rrr]^{\vee_{\L^nX}(\cdot)_n^{\deg
f}}_<(.5){\;}="h"&&&S^n_X\ar@/^20pt/[rd]^{ev}&
\\S^n\ar@/^20pt/[ru]^0_{\;}="b"\ar[r]_{\S^{n-1}\varepsilon}&
S^n\vee S^n\ar[u]_{(1,1)}^<(.1){\;}="a"\ar[r]_{\S^{n-1}(g\vee g)}&
S^n_X\vee S^n_X\ar[u]^{(1,1)}
\ar@/^20pt/[rrru]_<(.6){\;}="f"^<(.7){\;}="g"
\ar@/_30pt/[rrru]_<(.4){(\vee_{\L^nX}f,\vee_{\L^nX}(\cdot)_n^{\deg
f})}^<(.6){\;}="e"\ar[rrr]_{(\S^n f^*,\vee_{\L^nX}(\cdot)_n^{\deg
f})}
&&&S^n_X\ar[r]_{ev}&X\ar@{=>}"a";"b"^N\ar@{=>}"e";"f"^<(.1){(\vee_{\L^nX}\gamma,0^\vc)}
\ar@{=>}"g";"h"^{(\vee_{\L^nX}\delta,0^\vc)}}
\end{equation*}

Let us now pay special attention to the following subdiagram of (a)
\begin{equation*}\tag{b}
\xymatrix@R=100pt{&S^n\ar[r]^{\S^{n-1}g}&S^n_X\ar[rrr]^<(.1){\;}="c"^{\vee_{\L^nX}(\cdot)_n^{\deg
f}}_<(.5){\;}="h"&&&S^n_X&
\\S^n\ar@/^20pt/[ru]^0_{\;}="b"\ar[r]_{\S^{n-1}\varepsilon}&
S^n\vee S^n\ar[u]_{(1,1)}^<(.1){\;}="a"\ar[r]_{\S^{n-1}(g\vee g)}&
S^n_X\vee S^n_X\ar[u]^{(1,1)}
\ar@/^20pt/[rrru]_<(.6){\;}="f"^<(.7){\;}="g"_{(\vee_{\L^nX}(\cdot)_n^{\deg
f},\vee_{\L^nX}(\cdot)_n^{\deg f})} &&&&
\ar@{=>}"g";"h"^{(\vee_{\L^nX}\delta,0^\vc)}\ar@{=>}"a";"b"^N}
\end{equation*}
This is a composite track, termed (b), between $(n-1)$-fold
suspensions. The reduced Hopf invariant of (b) is trivial if $\overline{\hopf}(\delta)=0$. If
$0\neq\overline{\hopf}(\delta)\in\hat{\otimes}^2\Z=\Z/2$ then the reduced Hopf invariant of (b) is given by the following formula
where $(\pi_1g)_\abb(1)=\sum_{i=0}^k n_ia_i\in\Z[\L^nX]$ for some
$a_i\in \L^nX$ and $n_i\in\Z$,
\begin{equation*}\tag{c}
\overline{\hopf}\text{(b)}=\sum_{i=0}^kn_ia_i\hat{\otimes}a_i\in\hat{\otimes}^2\Z[\L^nX].
\end{equation*}
Here we have used the elementary properties of the Hopf invariant
for tracks described in \cite{2hg1} 3. By using again these properties
the reader can easily check that the following composite track has the same
reduced Hopf invariant as (b)
\begin{equation*}\tag{d}
\xymatrix@R=100pt{&&&S^n\ar[rrd]^{\S^{n-1}g}&&&&\\&S^n\vee
S^n\ar@/^20pt/[rru]|<(.13){\;\;\;((\cdot)^{\deg f}_n,(\cdot)^{\deg
f}_n)}_{\;}="j"\ar@/_20pt/[rru]^{\;}="i"|<(.7){((\cdot)^{\deg
f}_n,(\cdot)^{\deg
f}_n)}_{\;}="d"\ar[r]_{\S^{n-1}(g,g)}&S^n_X\ar[rrr]^<(.1){\;}="c"^{\vee_{\L^nX}(\cdot)_n^{\deg
f}}_<(.5){\;}="h"&&&S^n_X&
\\S^n\ar@/^60pt/[rrruu]^0_<(.2){\;}="b"\ar[r]_{\S^{n-1}\varepsilon}&
S^n\vee S^n\ar@{=}[u]^<(.1){\;}="a"\ar[r]_{\S^{n-1}(g\vee g)}&
S^n_X\vee S^n_X\ar[u]^{(1,1)} &&&&
\ar@{=>}"c";"d"_{(Q,Q)}\ar@{=>}"a";"b"^N\ar@{=>}"i";"j"^{(\delta,0^\vc)}}
\end{equation*}
Here $$Q\colon (\vee_{\L^nX}(\cdot)^{\deg f}_n)(\S^{n-1}g)\rr
(\S^{n-1}g)(\cdot)^{\deg f}_n$$ can be any 
track. Since (b) and (d) have the same reduced Hopf invariantthen (a) represents the same element in $\Pi_{n,1}X$
as 
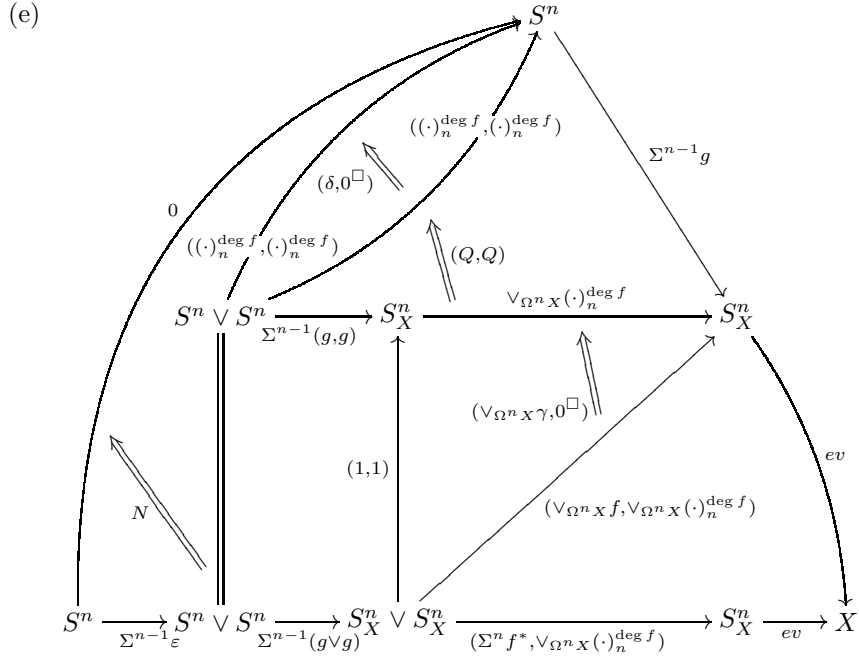
\begin{equation*}\tag{e}
\xymatrix@R=100pt{&&&S^n\ar[rrd]^{\S^{n-1}g}&&&&\\&S^n\vee
S^n\ar@/^20pt/[rru]|<(.13){\;\;\;((\cdot)^{\deg f}_n,(\cdot)^{\deg
f}_n)}_{\;}="j"\ar@/_20pt/[rru]^{\;}="i"|<(.7){((\cdot)^{\deg
f}_n,(\cdot)^{\deg
f}_n)}_{\;}="d"\ar[r]_{\S^{n-1}(g,g)}&S^n_X\ar[rrr]^<(.1){\;}="c"^{\vee_{\L^nX}(\cdot)_n^{\deg
f}}_<(.5){\;}="f"&&&S^n_X\ar@/^10pt/[rd]^{ev}&
\\S^n\ar@/^60pt/[rrruu]^0_<(.2){\;}="b"\ar[r]_{\S^{n-1}\varepsilon}&
S^n\vee S^n\ar@{=}[u]^<(.1){\;}="a"\ar[r]_{\S^{n-1}(g\vee g)}&
S^n_X\vee S^n_X\ar[u]^{(1,1)}
\ar[rrru]_<(.4){(\vee_{\L^nX}f,\vee_{\L^nX}(\cdot)_n^{\deg
f})}^<(.6){\;}="e"\ar[rrr]_{(\S^n f^*,\vee_{\L^nX}(\cdot)_n^{\deg
f})} &&&S^n_X\ar[r]_{ev}&X
\ar@{=>}"c";"d"_{(Q,Q)}\ar@{=>}"a";"b"^N\ar@{=>}"i";"j"^{(\delta,0^\vc)}\ar@{=>}"e";"f"^<(.1){(\vee_{\L^nX}\gamma,0^\vc)}}
\end{equation*}
The composite track (e) is the same as
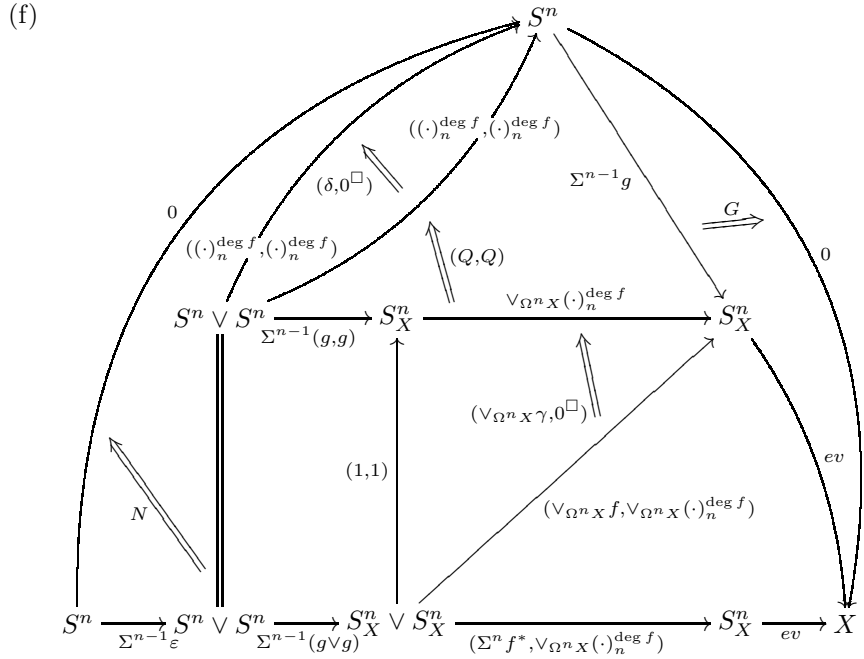
\begin{equation*}\tag{f}
\xymatrix@R=100pt{&&&S^n\ar@/^50pt/[rrrdd]^0_<(.4){\;}="l"\ar[rrd]^<(.7){\;}="k"_{\S^{n-1}g}&&&&\\&S^n\vee
S^n\ar@/^20pt/[rru]|<(.13){\;\;\;((\cdot)^{\deg f}_n,(\cdot)^{\deg
f}_n)}_{\;}="j"\ar@/_20pt/[rru]^{\;}="i"|<(.7){((\cdot)^{\deg
f}_n,(\cdot)^{\deg
f}_n)}_{\;}="d"\ar[r]_{\S^{n-1}(g,g)}&S^n_X\ar[rrr]^<(.1){\;}="c"^{\vee_{\L^nX}(\cdot)_n^{\deg
f}}_<(.5){\;}="f"&&&S^n_X\ar@/^10pt/[rd]^{ev}&
\\S^n\ar@/^60pt/[rrruu]^0_<(.2){\;}="b"\ar[r]_{\S^{n-1}\varepsilon}&
S^n\vee S^n\ar@{=}[u]^<(.1){\;}="a"\ar[r]_{\S^{n-1}(g\vee g)}&
S^n_X\vee S^n_X\ar[u]^{(1,1)}
\ar[rrru]_<(.4){(\vee_{\L^nX}f,\vee_{\L^nX}(\cdot)_n^{\deg
f})}^<(.6){\;}="e"\ar[rrr]_{(\S^n f^*,\vee_{\L^nX}(\cdot)_n^{\deg
f})} &&&S^n_X\ar[r]_{ev}&X
\ar@{=>}"c";"d"_{(Q,Q)}\ar@{=>}"a";"b"^N\ar@{=>}"i";"j"^{(\delta,0^\vc)}\ar@{=>}"e";"f"^<(.1){(\vee_{\L^nX}\gamma,0^\vc)}
\ar@{=>}"k";"l"^G}
\end{equation*}
And (f) is the same as
\begin{equation*}\tag{g}
\xymatrix@R=100pt{&&&S^n\ar@/^50pt/[rrrdd]^0_<(.4){\;}="l"\ar[rrd]^<(.7){\;}="k"_{\S^{n-1}g}&&&&\\&S^n\vee
S^n\ar@/_20pt/[rru]^{\;}="i"|<(.7){((\cdot)^{\deg f}_n,(\cdot)^{\deg
f}_n)}_{\;}="d"\ar[r]_{\S^{n-1}(g,g)}&S^n_X\ar[rrr]^<(.1){\;}="c"^{\vee_{\L^nX}(\cdot)_n^{\deg
f}}_<(.5){\;}="f"&&&S^n_X\ar@/^10pt/[rd]^{ev}&
\\S^n\ar[r]_{\S^{n-1}\varepsilon}&
S^n\vee S^n\ar@{=}[u]^<(.1){\;}="a"\ar[r]_{\S^{n-1}(g\vee g)}&
S^n_X\vee S^n_X\ar[u]^{(1,1)}
\ar[rrru]_<(.4){(\vee_{\L^nX}f,\vee_{\L^nX}(\cdot)_n^{\deg
f})}^<(.6){\;}="e"\ar[rrr]_{(\S^n f^*,\vee_{\L^nX}(\cdot)_n^{\deg
f})} &&&S^n_X\ar[r]_{ev}&X
\ar@{=>}"c";"d"_{(Q,Q)}\ar@{=>}"e";"f"^<(.1){(\vee_{\L^nX}\gamma,0^\vc)}
\ar@{=>}"k";"l"^G}
\end{equation*}
Obviously (g) coincides with
\begin{equation*}\tag{h}
\xymatrix@R=100pt{&&&S^n\ar@/^50pt/[rrrdd]^0_<(.4){\;}="l"\ar[rrd]^<(.7){\;}="k"_{\S^{n-1}g}&&&&\\&
S^n\ar[rru]^{\;}="i"|<(.7){(\cdot)^{\deg
f}_n}_{\;}="d"\ar[r]^{\S^{n-1}g}&S^n_X\ar[rrr]^<(.1){\;}="c"^{\vee_{\L^nX}(\cdot)_n^{\deg
f}}_<(.5){\;}="f"&&&S^n_X\ar@/^10pt/[rd]^{ev}&
\\S^n\ar@/^20pt/[ru]^0_{\;}="b"\ar[r]_{\S^{n-1}\varepsilon}&
S^n\vee S^n\ar[u]_{(1,1)}^<(.1){\;}="a"\ar[r]_{\S^{n-1}(g\vee g)}&
S^n_X\vee S^n_X\ar[u]^{(1,1)}
\ar[rrru]_<(.4){(\vee_{\L^nX}f,\vee_{\L^nX}(\cdot)_n^{\deg
f})}^<(.6){\;}="e"\ar[rrr]_{(\S^n f^*,\vee_{\L^nX}(\cdot)_n^{\deg
f})} &&&S^n_X\ar[r]_{ev}&X
\ar@{=>}"c";"d"_Q\ar@{=>}"e";"f"^<(.1){(\vee_{\L^nX}\gamma,0^\vc)}
\ar@{=>}"k";"l"^G\ar@{=>}"a";"b"^N}
\end{equation*}
And (h) is the same as
\begin{equation*}\tag{i}
\xymatrix@R=100pt{&
S^n\ar[r]^{\S^{n-1}g}&S^n_X\ar[rrr]^<(.1){\;}="c"^{\vee_{\L^nX}(\cdot)_n^{\deg
f}}_<(.5){\;}="f"&&&S^n_X\ar@/^10pt/[rd]^{ev}&
\\S^n\ar@/^20pt/[ru]^0_{\;}="b"\ar[r]_{\S^{n-1}\varepsilon}&
S^n\vee S^n\ar[u]_{(1,1)}^<(.1){\;}="a"\ar[r]_{\S^{n-1}(g\vee g)}&
S^n_X\vee S^n_X\ar[u]^{(1,1)}
\ar[rrru]_<(.4){(\vee_{\L^nX}f,\vee_{\L^nX}(\cdot)_n^{\deg
f})}^<(.6){\;}="e"\ar[rrr]_{(\S^n f^*,\vee_{\L^nX}(\cdot)_n^{\deg
f})} &&&S^n_X\ar[r]_{ev}&X
\ar@{=>}"e";"f"^<(.1){(\vee_{\L^nX}\gamma,0^\vc)} \ar@{=>}"a";"b"^N}
\end{equation*}
Notice that this last composite track (i) represents
$\grupo{\partial[g,G],\gamma}$, hence we are done.
\end{proof}

The next corollary follows from the two previous propositions.

\begin{cor}\label{mona}
For any pointed space $X$ and $n\geq 2$ the monoid
$\pi_{0,0}\endo_*(S^n)$ acts on the right of the quadratic pair
module $\Pi_{n,*}X$. This action is natural in $X$. 
\end{cor}


Now Theorem \ref{laa} is a consequence of the next result.

\begin{prop}
The action of $\pi_{0,0}\endo_*(S^n)$ on the right of $\Pi_{n,*}X$
given by Corollary \ref{mona} extends to an action of the whole
monoid-groupoid $\pi_{0,*}\endo_*(S^n)$, $n\geq 2$.
\end{prop}

\begin{proof}
A morphism $H$ in $\pi_{0,*}\endo_*(S^n)$ is a track $H\colon
f\rr g$ between maps $f,g\colon S^n\r S^n$, in particular $\deg
f=\deg g=k\in\Z$. In order to define a track $$H^*\colon f^*\rr
g^*$$ between the quadratic pair module morphisms
$$f^*,g^*\colon\Pi_{n,*}X\To\Pi_{n,*}X$$
we choose tracks in $\tilde{\pi}_{0,1}\endo_*(S^n)$
$$\alpha\colon
f\rr(\cdot)_n^k,$$
$$\beta\colon
g\rr(\cdot)^k_n,$$ such that
$$H=\beta^\vi\vc\alpha.$$
By Proposition \ref{nores} the maps
$$\grupo{-,\alpha},\grupo{-,\beta}\colon\Pi_{n,0}X\To\Pi_{n,1}X$$
are tracks
$$\grupo{-,\alpha}\colon f^*\rr k^*,$$
$$\grupo{-,\beta}\colon g^*\rr k^*,$$
in the category $\C{wqpm}$, therefore we can define
$H^*$ as the vertical composition
$$H^*=\grupo{-,\beta}^\vi\vc\grupo{-,\alpha},$$
i. e. $H^*$ is the map
$$H^*\colon\Pi_{n,0}X\To\Pi_{n,1}X$$
defined by
$$H^*(x)=\grupo{x,\alpha}-\grupo{x,\beta}.$$
By the proof of Proposition \ref{nores} and by \cite{2hg1} 7.5 (3) the element $H(x)$ coincides with
$r(ev(\vee_{\L^nX}H)(\S^{n-1}\tilde{x}))$ for $\tilde{x}\colon S^1\r \vee_{\L^nX}S^1$ any map with
$(\pi_1\tilde{x})_\ni(1)=x$
in the sense of \cite{2hg1} 7.4.
The reader can now use the properties of the bracket (\ref{bra})
described in Proposition \ref{nores} together with \cite{2hg1} 7.5 (3) to check that this yields a
monoid-groupoid action.
\end{proof}



Later we will consider the quotient monoid $\bar{\pi}_{0,1}\endo_*(S^2)$ of $\tilde{\pi}_{0,1}\endo_*(S^2)$
defined as follows: two elements $\gamma\colon f\rr(\cdot)_2^{\deg f}$, $\bar{\gamma}\colon g\rr(\cdot)_2^{\deg g}$
in $\tilde{\pi}_{0,1}\endo_*(S^2)$ represent the same element in $\bar{\pi}_{0,1}\endo_*(S^2)$ provided $\deg f=\deg g$ and
$$0=\overline{\hopf}(\bar{\gamma}\vc\gamma^\vi)\in\hat{\otimes}^2\Z=\Z/2.$$

\begin{prop}\label{fatora}
The bracket
operation (\ref{bra}) factors for $n=2$ through the natural projection
$\tilde{\pi}_{0,1}\endo_*(S^2)\twoheadrightarrow\bar{\pi}_{0,1}\endo_*(S^2)$. 
$$\grupo{-,-}\colon\Pi_{2,0}\times\bar{\pi}_{0,1}\endo_*(S^2)\To\Pi_{2,1}X.$$
\end{prop}

\begin{proof}
Two tracks $\gamma$ and $\bar{\gamma}$ in $\tilde{\pi}_{0,1}\endo_*(S^2)$ represent the same element in 
$\bar{\pi}_{0,1}\endo_*(S^2)$ if and only if $\bar{\gamma}=\delta\vc\gamma$ for some 
$\delta\colon(\cdot)^k_2\rr(\cdot)^k_2$ with $\overline{\hopf}(\delta)=0$, so we only need to check that 
$\grupo{x,\gamma}=\grupo{x,\delta\vc\gamma}$. The element $\grupo{x,\delta\vc\gamma}$ is represented by 
diagram (a) in the proof of Proposition \ref{oes} where we assume that $g$ is a map with $(\pi_1g)_\ni(1)=x$. 
As we mention in that proof diagram (b) is a nil-track in these circumstances, therefore we can drop $\delta$ 
from (a) and still obtain the same element in $\Pi_{2,1}X$. But if we drop $\delta$ we obtain $\grupo{x,\gamma}$, hence we are done.
\end{proof}

\section{The symmetric action on secondary homotopy
groups}\label{the}


The permutation of coordinates in
$S^n=S^1\wedge\cdots\wedge S^1$
induces a left action of the symmetric group $\sym{n}$ on the $n$-sphere $S^n$. This action induces a monoid
inclusion
\begin{equation}\label{sima}
\sym{n}\subset\pi_{0,0}\endo_*(S^n).
\end{equation}

We define the \emph{symmetric track group} for $n\geq 3$
$$\symt{n}\subset\tilde{\pi}_{0,1}\endo_*(S^n)$$
as the submonoid of tracks of the from $$\alpha\colon
\sigma\rr(\cdot)_n^{\sign(\sigma)},$$ where $\sigma\in\sym{n}$ and
$\sign(\sigma)\in\set{\pm1}$ is the sign of the permutation. Compare
the notation in (\ref{nse}) and (\ref{nseya}).

The submonoid defined as above for $n=2$ will be called the \emph{extended symmetric track group}
$$\symtex{2}\subset\tilde{\pi}_{0,1}\endo_*(S^2).$$
For $n=2$ the symmetric track group $\symt{2}$ is the image of $\symtex{2}$ by the natural projection 
$\tilde{\pi}_{0,1}\endo_*(S^2)\twoheadrightarrow\bar{\pi}_{0,1}\endo_*(S^2)$ in Proposition \ref{fatora}.
\begin{equation}\label{sima2}
\symt{2}\subset\bar{\pi}_{0,1}\endo_*(S^2).
\end{equation}

\begin{prop}\label{ps}
The symmetric track group is indeed a group. Moreover, it fits into
a central extension, $n\geq 2$,
$$\Z/2\hookrightarrow\symt{n}\st{\delta}\twoheadrightarrow\sym{n}$$
with $\delta(\alpha)=\sigma$, which splits if and only if $n=2$ or $3$.
\end{prop}

This proposition follows from Corollary \ref{yi} and Remarks \ref{yak} and \ref{n2} below.

For $n=0$ and $n=1$ we define $\sym{n}$ to be the trivial group, and $\symt{n}$ the trivial sign group.
Then the symmetric track group $\symt{n}$ is a sign group $(n\geq
0)$
$$\set{\pm1}\hookrightarrow\symt{n}\st{\delta}\twoheadrightarrow\sym{n}\st{\sign}\To\set{\pm1}$$
as in Definition \ref{signg}.

\begin{thm}
Let $X$ be a pointed space. For $n\geq 0$ the symmetric group
$\sym{n}$ acts naturally on the right of the additive secondary homotopy
group $\Pi_{n,*}X$ in the category $\C{qpm}$ of quadratic pair
modules. Moreover, the restriction
$$\grupo{-,-}\colon\Pi_{n,0}X\times\symt{n}\To\Pi_{n,1}X$$
of the bracket defined in (\ref{bra}) if $n\geq 3$ and in Proposition \ref{fatora} if $n=2$ yields a natural right action
of the sign group $\symt{n}$ on $\Pi_{n,*}X$ in the sense of Definition \ref{signg}.
\end{thm}

The action of $\sym{n}$ is given by Corollary \ref{mona} and the
inclusion (\ref{sima}) if $n\geq 3$ or (\ref{sima2}) if $n=2$. The rest of the statement follows from
 Proposition \ref{nores}. The cases $n=0,1$ are trivial consequences of Remark \ref{pursu}.

\section{The structure of the symmetric track groups}\label{krio}

In this section we construct a positive pin
representation for the symmetric track group $\symt{n}$. By using this representation we obtain a finite 
presentation of $\symt{n}$.

The action of $\sym{n}$ on $S^n$ can be extended to a well-known
action of the orthogonal group $O(n)$ which we now recall. Let $[-1,1]^n\subset\Real^n$
be the hypercube centered in the origin whose vertices have all
coordinates in $\set{\pm1}$, $D^n\subset\Real^n$ the Euclidean unit
ball and $S^{n-1}$ its boundary. There is a homeomorphism
$\phi\colon[-1,1]^n\r D^n$
fixing the origin defined as follows
$$\phi(\ul{x})=\frac{\max\limits_{1\leq i\leq n}|x_i|}{\norm{\ul{x}}}\;\ul{x}.$$
Here $\ul{x}\in[-1,1]^n$ is an arbitrary non-trivial vector in the
hypercube and $\norm{\cdot}$ is the Euclidean norm. This
homeomorphism projects the hypercube onto the ball from the origin.
There is also a map collapsing the boundary
$$\varrho\colon[-1,1]^n\To S^1\wedge\st{n}\cdots\wedge S^1=S^n,$$
$$\varrho(x_1,\dots,x_n)=(\exp(i\pi(1+x_1)),\dots,\exp(i\pi(1+x_n))).$$
The composite $$\varrho\phi^{-1}\colon D^n\To S^n$$
induces a homeomorphism
$$D^n/S^{n-1}\cong S^1\wedge\st{n}\cdots\wedge S^1=S^n$$
that we fix.

The orthogonal group $O(n)$ acts on the left of the unit ball $D^n$.
This action induces an action of $O(n)$ on the quotient space
$S^n=D^n/S^{n-1}$ preserving the base-point. 
The interchange of coordinates action of the symmetric group
$\sym{n}$ on $\Real^n$ preserves the Euclidean scalar product, and
therefore induces a homomorphism 
\begin{equation}\label{lai}
i\colon\sym{n}\hookrightarrow
O(n).
\end{equation} The pull-back of the action of $O(n)$ along this
homomorphism is the action of $\sym{n}$ on $S^n$ given by the smash
product decomposition of $S^n$.

\begin{rem}\label{esJ}
The action of $O(n)$ on $S^n$ defines an inclusion
$O(n)\subset\endo_*(S^n)$. The induced homomorphism on $\pi_1$ is
the Whitehead-Hopf $J$-homomorphism
\begin{equation}\label{J}
J\colon\pi_1O(n)\cong\pi_1\endo_*(S^n)=\pi_{n+1}S^n
\end{equation}
which is known to be an isomorphism for $n\geq 2$. Let $\pi_{0,*}O(n)$ be the fundamental groupoid of the Lie
group $O(n)$. Then, considering elements $A,B\in O(n)$ as pointed maps
$$A,B\colon S^n\To S^n$$
the isomorphism in (\ref{J}) allows to identify all morphisms $\gamma\colon A\r B$ in $\pi_{0,*}O(n)$ 
with all tracks
$$\gamma\colon A\rr B$$
in $\pi_{0,1}\endo_*(S^n)$. 
Let $\id_n\in O(n)$ be the identity
matrix. The order $2$ matrix
$$\left(
    \begin{array}{cc}
      \id_{n-1} & 0 \\
      0 & -1 \\
    \end{array}
  \right)\in O(n)
$$
will be denoted by $\id_{n-1}\oplus(-1)$. By using the action of
$O(n)$ on $S^n$ we have by the notation in (\ref{nseya}) that
$$\id_{n-1}\oplus(-1)=(\cdot)^{-1}_n\colon S^n\To S^n.$$
Obviously $\id_n=(\cdot)^ 1_n=1_{S^n}\colon S^n\r S^n$.
\end{rem}

The topological group structure of $O(n)$ induces an internal group
structure on the fundamental groupoid $\pi_{0,*}O(n)$ in the
category of groupoids. In particular the set $\pi_{0,1}O(n)$ of
morphisms in $\pi_{0,*}O(n)$ forms a group. We define the subgroup
$$\widetilde{O}(n)\subset\pi_{0,1}O(n)$$ consisting of all the morphisms
with target $\id_n$ or $\id_{n-1}\oplus (-1)$. 
By Remark \ref{esJ} the symmetric track group is the subgroup 
$$\symt{n}\subset\widetilde{O}(n)$$ of morphisms with source in the image of $i$ in (\ref{lai}), $n\geq 3$.
The subgroup $\widetilde{O}(n)$ is
embedded in an extension
\begin{equation}\label{pin}
\Z/2\hookrightarrow  \widetilde{O}(n)\st{q}\twoheadrightarrow
O(n),\;\;\; n\geq 3.
\end{equation}
The projection $q$ sends a morphism in
$\widetilde{O}(n)\subset\pi_{0,1}O(n)$ to the source, and the kernel
is clearly
$\pi_1O(n)=\Z/2$ for $n\geq 3$. The case $n=2$ will be considered in Remark \ref{n2} below.

There is also a well-known extension
\begin{equation}\label{pinext}
\Z/2\hookrightarrow Pin_+(n)\st{\rho}\twoheadrightarrow
O(n)
\end{equation}
given by the positive pin group. Let us recall the definition of
this extension.

\begin{defn}\label{clif}
The positive Clifford algebra $C_+(n)$ is the unital $\Real$-algebra
generated by $e_i$, $1\leq i\leq n$, with relations
\begin{enumerate}
\item $e_i^2=1$ for $1\leq i \leq n$,
\item $e_ie_j=-e_je_i$ for $1\leq i<j \leq n$.
\end{enumerate}
Clifford algebras are defined for arbitrary quadratic forms on
finite-dimensional vector spaces, see for instance \cite{rclg} 6.1.
The Clifford algebra defined above corresponds to the quadratic form
of the standard positive-definite scalar product in $\Real^n$. We
identify the sphere $S^{n-1}$ with the vectors of Euclidean norm
$1$ in the vector subspace $\Real^n\subset C_+(n)$ spanned by the
generators $e_i$. The vectors in $S^{n-1}$ are units in $C_+(n)$.
Indeed for any $v\in S^{n-1}$ the square $v^2=1$ is the unit element
in $C_+(n)$, so that $v^{-1}=v$. The group $Pin_+(n)$ is the
subgroup of units in $C_+(n)$ generated by $S^{n-1}$. Any $x\in
Pin_+(n)$ defines an automorphism of $\Real^n\subset C_+(n)$ given
by conjugation in $C_+(n)$ as follows
$$\Real^n\To\Real^n\colon w\mapsto -xwx^{-1}.$$
If $x\in S^{n-1}$ then this automorphism is the reflection along the
hyperplane orthogonal to the unit vector $x$. This endomorphism
always preserves the scalar product, therefore this
defines a homomorphism
$$\rho\colon Pin_+(n)\twoheadrightarrow O(n).$$
This homomorphism is surjective since all elements in $O(n)$ are
products of $\leq n$ reflections. It is easy to see that the kernel
of $\rho$ is $\Z/2$ generated by $-1\in C_+(n)$. This is the
extension in (\ref{pinext}).

The Clifford algebra $C_+(n)$ has dimension $2^n$. A basis is given
by the elements $$e_{i_1}\cdots e_{i_k},\;\;1\leq i_1<\cdots<i_k\leq
n.$$ We give $C_+(n)$ the topology induced by the Euclidean norm
associated to this basis. The positive pin group inherits a topology
turning (\ref{pinext}) into a Lie group extension.
\end{defn}

\begin{prop}\label{mp}
The extension (\ref{pin}) is isomorphic to (\ref{pinext}).
\end{prop}

\begin{proof}

Since $Pin_+(n)$ is a topological group $\pi_{0,*}Pin_+(n)$ is a
group object in the category of groupoids. We define
$$\widetilde{Pin}_+(n)\subset\pi_{0,1}Pin_+(n)$$
to be the subgroup given by morphisms $x\r y$ in $\pi_{0,*}Pin_+(n)$
with target $1$ or $e_n$. This is well defined since
$\set{1,e_n}\subset Pin_+(n)$ is a subgroup. This observation is indeed the key step of the proof, and it shows
for example why the negative pin group does not occur as (\ref{pin}). Moreover,
$\pi_{0,*}\rho$ induces a homomorphism
\begin{equation*}\tag{a}
\widetilde{Pin}_+(n)\To \widetilde{O}(n).
\end{equation*}

It is well-known that $Pin_+(n)$ has two components. The two
components are separated by the function
$$Pin_+(n)\st{\rho}\twoheadrightarrow O(n)\st{\det}\twoheadrightarrow \set{\pm1}.$$
In particular $1$ and $e_n$ lie in different components, hence the
homomorphism
\begin{equation*}\tag{b}
\widetilde{Pin}_+(n)\To Pin_+(n)
\end{equation*}
is surjective. Moreover, it is injective since the two components of
$Pin_+(n)$ are known to be simply connected, therefore (b) is an
isomorphism. The inverse
\begin{equation*}\tag{c}
Pin_+(n)\To\widetilde{Pin}_+(n)
\end{equation*}
sends an element $x\in Pin_+(n)$ to the image by $\pi_{0,*}\rho$ of
the unique morphism $x\r y$ in $\pi_{0,*}Pin_+(n)$ with $y=e_n$
provided $\det\rho(x)=-1$ or $y=1$ otherwise.

Obviously the composite of (c) and (a) is compatible with the
projections onto $O(n)$ in (\ref{pin}) and (\ref{pinext}), so we
only need to check that the composite of (c) and (a) induces an
isomorphism between the kernels. The kernel of $\rho$ is $-1$. A
path $\gamma\colon[0,1]\r Pin_+(n)$ from $-1$ to $1$ is defined by
$$\gamma(t)=(-\cos(t\pi)e_2+\sin(t\pi)e_1)e_2=-\cos(t\pi)+\sin(t\pi)e_1e_2.$$
Now it is an easy exercise  to check that
$\rho\gamma\colon [0,1]\r O(n)$ is a generator of $\pi_1O(n)$, and
hence we are done.

\end{proof}

\begin{rem}
We recall that the extension (\ref{pinext}), and therefore
(\ref{pin}), represents the second Stiefel-Whitney class $w_2\in
H^2(BO(n),\Z/2)$, compare \cite{t4mffg} page 21.
\end{rem}

By definition of (\ref{pin}) and Proposition \ref{mp} we obtain the
following corollary.

\begin{cor}\label{yi}
For $n\geq 3$ the symmetric track group $\symt{n}$ is the pull back
of the central extension for the positive pin group $Pin_+(n)$ in
(\ref{pinext}) along the inclusion $i\colon\sym{n}\subset O(n)$, in
particular there is a central extension
$$\Z/2\hookrightarrow \symt{n}\st{\delta}\twoheadrightarrow\sym{n}$$
classified by the pull-back of the second Stiefel-Whitney class
$i^*w_2\in H^2(\sym{n},\Z/2)$. 
\end{cor}

\begin{rem}\label{yak}
The low-dimensional mod $2$ cohomology groups of symmetric groups
$\sym{n}$ are as follows, $n\geq 3$,
$$H^1(\sym{n},\Z/2)=\left\{
                      \begin{array}{ll}
                        \Z/2\, \chi\oplus\Z/2\, i^*w_1, & \hbox{for $n=3$;} \\
                        &\\
                        \Z/2\, i^*w_1, & \hbox{for $n>3$;}
                                              \end{array}
                    \right.
$$
$$H^2(\sym{n},\Z/2)=\left\{
                      \begin{array}{ll}
                        \Z/2\, i^*w_1^2, & \hbox{for $n=3$;} \\
                        &\\
                        \Z/2\, i^*w_1^2\oplus\Z/2\, i^*w_2, & \hbox{for $n>3$.}
                                              \end{array}
                    \right.
$$
Here we write $w_j\in H^j(BO(n),\Z/2)$ for the $j^\text{th}$
Stiefel-Whitney class, $j=1,2$. The pull-back $i^*w_1$ corresponds
to the sign homomorphism
$$i^*w_1=\sign\colon\sym{n}\To \set{\pm1}\cong\Z/2,$$
The pull-back of the second Stiefel-Whitney class is trivial for
$n=3$, therefore $\symt{3}$ is a split extension of $\sym{3}$ by
$\Z/2$, and $\chi\colon\symt{3}\twoheadrightarrow\Z/2$ is a
retraction.
\end{rem}

The following structure theorem follows from Corollary
\ref{yi}.

\begin{thm}\label{ultimate}
The symmetric track group $\symt{n}$ is the subgroup of $Pin_+(n)$
formed by the units $x\in C_+(n)$ such that for any $1\leq i\leq n$
there exists $1\leq\sigma(i)\leq n$ with
$-xe_ix^{-1}=e_{\sigma(i)}$. The boundary homomorphism
$\delta\colon\symt{n}\twoheadrightarrow\sym{n}$ sends $x$ above to
the permutation $\delta(x)=\sigma$. The group $\symt{n}$ has a
presentation given by generators $\omega$, $t_i$, $1\leq i\leq n-1$,
and relations
\begin{eqnarray*}
t_1^2&=&1\text{ for }1\leq i\leq n-1,\\
(t_it_{i+1})^3&=&1\text{ for }1\leq i\leq n-2,\\
\omega^2&=&1,\\
t_i\omega&=&\omega t_i\text{ for }1\leq i\leq n-1,\\
t_it_j&=&\omega t_jt_i\text{ for }1\leq i<j-1\leq n-1;
\end{eqnarray*}
with $\omega\mapsto -1$ and $t_i\mapsto \frac{1}{\sqrt{2}}(e_i-e_{i+1})$. In particular
$\delta(\omega)=0$ and $\delta(t_i)=(i\; i+1)$.
\end{thm}

This is a group considered by Schur in \cite{dsag} and by Serre in
\cite{iwf}.

\begin{rem}\label{n2}
In case $n=2$ we have $O(2)=\set{\pm1}\ltimes S^1$ with $\set{\pm1}$ acting on $S^1$ exponentially,
$\widetilde{O}(2)=\set{\pm1}\ltimes\Real$ with $\set{\pm1}$ acting on $\Real$ multiplicatively, and the
projection $q\colon\widetilde{O}(2)\twoheadrightarrow O(2)$ defined as in (\ref{pin}) is the identity in
$\set{\pm1}$ and the exponential map in the second coordinate $\Real\twoheadrightarrow S^1\colon x\mapsto \exp
(2\pi i x)$. In particular we have an abelian extension
\begin{equation*}
\Z\hookrightarrow\widetilde{O}(2)\st{q}\twoheadrightarrow O(2).
\end{equation*}
The induced action of $O(2)$ on $\Z$ is given by the determinant $\det\colon O(2)\twoheadrightarrow\set{\pm1}$. 
By Remark \ref{esJ} the extended
symmetric track group $\symtex{2}$ is the pull-back of $i\colon\sym{2}\subset O(2)$ along $q$, therefore we have an
abelian extension
\begin{equation}\label{pina}
\Z\hookrightarrow\symtex{2}\st{q}\twoheadrightarrow \sym{2},
\end{equation}
where $\symt{2}$ acts on $\Z$ by the unique isomorphism $\symt{2}\cong\set{\pm1}$. Now the symmetric track group
$\symt{2}$ can be identified with the push-forward of the extension (\ref{pina}) along the natural projection
$\Z\twoheadrightarrow\Z/2$, therefore we get a central extension
\begin{equation}\label{pinako}
\Z/2\hookrightarrow\symt{2}\st{q}\twoheadrightarrow \sym{2}.
\end{equation}
The cohomology group
$H^2(\symt{2},\Z)=0$ is trivial, so (\ref{pina}) is a splitting extension. Moreover (\ref{pinako}) is also
splitting since it is the push-forward of (\ref{pina}).
\end{rem}

\section{An application to the cup-one-product}

Let $n\geq m> 1$ be even integers. The cup-one-product operation
$$\pi_nS^m\To\pi_{2n+1}S^{2m}\colon\alpha\mapsto\alpha\cup_1\alpha$$
is defined in the following way, compare \cite{c1Tb} 2.2.1. Let $k$ be any positive integer and let 
$\tau_k\in\sym{2k}$ be the permutation exchanging the first and
the second block of $k$ elements in $\set{1,\dots,2k}$. If $k$ is even then $\sign \tau_k=1$. 
We choose for any even integer $k>1$ a track
$\hat{\tau}_k\colon\tau_k \rr1_{S^{2k}}$ in $\symt{2k}$. Consider the following diagram in the track category
$\C{Top}^*$ of pointed spaces where 
$a\colon S^n\r S^m$ represents $\alpha$.
\begin{equation}\label{cup1}
\xymatrix{S^{2n}\ar[r]^{a\wedge a}\ar[d]^{{\tau}_n}_{\;}="b"\ar@/_30pt/[d]_{1_{S^{2n}}}^{\;}="a"&
S^{2m}\ar[d]_{{\tau}_m}^{\;}="c"\ar@/^30pt/[d]^{1_{S^{2m}}}_{\;}="d"\\
S^{2n}\ar[r]^{a\wedge a}&S^{2m}\ar@{=>}"a";"b"^{\hat{\tau}_n^\vi}\ar@{=>}"c";"d"^{\hat{\tau}_m}}
\end{equation}
By pasting this diagram we obtain a self-track of $a\wedge a$
\begin{equation}\label{cup2}
(\hat{\tau}_m(a\wedge a))\vc((a\wedge a)\hat{\tau}_n^\vi)\colon a\wedge a\rr a\wedge a.
\end{equation}
The set of self-tracks $a\wedge a\rr a\wedge a$ is the automorphism group of the map $a\wedge a$ in the track
category $\C{Top}^*$. 
The element $\alpha\cup_1\alpha\in\pi_{2n+1}S^{2m}$ is given by the track (\ref{cup2})
via the well-known Barcus-Barratt-Rutter isomorphism
$$\aut(a\wedge a)\cong\pi_{2n+1}S^{2m},$$
see \cite{hcefm}, \cite{hcmifs} and also \cite{ch4c} VI.3.12 and \cite{slort} for further details.

The following theorem generalizes \cite{K1p} 6.5.

\begin{thm}\label{cupthm}
The formula 
$$2(\alpha\cup_1\alpha)=\frac{n+m}{2}(\alpha\wedge\alpha)(\S^{2(n-1)}\eta)$$
holds, 
where $\eta\colon S^3\r S^2$ is the Hopf map.
\end{thm}

The proof of Theorem \ref{cupthm} is based on the following lemma.

\begin{lem}\label{A}
The following formula holds in $\symt{2k}$
$$\hat{\tau}^2_k=\omega^{\binom{k}{2}}.$$
\end{lem}

\begin{proof}
Here we use the representation of $\symt{2k}$ in $Pin_+(2k)$ given in
Theorem \ref{ultimate} and the relations (1) and (2) in the definition
of the Clifford algebra $C_+(2k)$, see Definition \ref{clif}.

The permutation $\tau_k$ can be expressed as a product of
transpositions as follows
$$\tau_k=(1\;\; k)(2\;\; k+1)\cdots(k-1\;\; 2k-1)(k\;\; 2k).$$

The element $\frac{1}{\sqrt{2}}(e_i- e_{i+k})\in S^{2k-1}\subset
Pin_+(2k)$ acts on $\Real^{2k}$ (with coordinates $e_i$, $1\leq i\leq
2k$) by reflection along the plane orthogonal to
$\frac{1}{\sqrt{2}}(e_i- e_{i+k})$, see Definition \ref{clif}. This
pane is $e_i=e_{i+k}$, therefore the action of
$\frac{1}{\sqrt{2}}(e_i- e_{i+k})$ on $\Real^{2k}$ interchanges the
coordinates in $e_i$ and $e_{i+k}$ and preserves all the other ones.

Now by Theorem \ref{ultimate} $\frac{1}{\sqrt{2}}(e_i- e_{i+k})$
lies in the positive pin representation of $\symt{2k}$ and
$\delta(\frac{1}{\sqrt{2}}(e_i- e_{i+k}))=(i\;\; i+k)$, so 
$$\hat{\tau}_k=\pm\frac{1}{2^{\frac{k}{2}}}(e_1-e_{k+1})(e_2-e_{k+2})\cdots(e_{k-1}-e_{2k-1})(e_k-e_{2k}).$$

The following equalities hold in the Clifford algebra $C_+(2k)$, see
the defining relations in Definition \ref{clif}, $i\neq j$,
\begin{eqnarray*}
(e_i-e_{i+k})^2&=&e_i^2-e_ie_{i+k}-e_{i+k}e_i+e_{i+k}^2\\
&=&1-e_ie_{i+k}+e_ie_{i+k}+1\\
&=&2,\\
(e_i-e_{i+k})(e_j-e_{j+k})&=&e_ie_j-e_ie_{j+k}-e_{i+k}e_j+e_{i+k}e_{j+k}\\
&=&-e_je_i+e_{j+k}e_i+e_je_{i+k}-e_{j+k}e_{i+k}\\
&=&-(e_j-e_{j+k})(e_i-e_{i+k}).
\end{eqnarray*}
Hence we observe that
\begin{eqnarray*}
\hat{\tau}_k^2&=&\frac{1}{2^{\frac{k}{2}}}\frac{1}{2^{\frac{k}{2}}}(-1)^{k-1}2(-1)^{k-2}2
\cdots (-1)^12(-1)^02\\
&=&(-1)^{k-1}(-1)^{k-2}\cdots(-1)^1(-1)^0\\
&=&(-1)^{\binom{k}{2}}.
\end{eqnarray*}
The proof is now finished.
\end{proof}

\begin{proof}[Proof of Theorem \ref{cupthm}]
The element $2(\alpha\cup_1\alpha)$ corresponds to the pasting of the following diagram 
$$\xymatrix{S^{2n}\ar[r]^{a\wedge a}\ar[d]^{{\tau}_n}_{\;}="f"\ar@/_30pt/[d]_{1_{S^{2n}}}^{\;}="e"&
S^{2m}\ar[d]_{{\tau}_m}^{\;}="g"\ar@/^30pt/[d]^{1_{S^{2m}}}_{\;}="h"\\
S^{2n}\ar[r]^{a\wedge a}\ar[d]^{{\tau}_n}_{\;}="b"\ar@/_30pt/[d]_{1_{S^{2n}}}^{\;}="a"&
S^{2m}\ar[d]_{{\tau}_m}^{\;}="c"\ar@/^30pt/[d]^{1_{S^{2m}}}_{\;}="d"\\
S^{2n}\ar[r]^{a\wedge a}&S^{2m}\ar@{=>}"a";"b"^{\hat{\tau}_n^\vi}\ar@{=>}"c";"d"^{\hat{\tau}_m}
\ar@{=>}"e";"f"^{\hat{\tau}_n^\vi}\ar@{=>}"g";"h"^{\hat{\tau}_m}}$$
By Lemma \ref{A} and using that $n$ and $m$ are even this composite track coincides with
$$\xymatrix{S^{2n}\ar@/^20pt/[rr]^{1_{S^{2n}}}_{\;}="a"\ar@/_20pt/[rr]_{1_{S^{2n}}}^{\;}="b"
&&S^{2n}\ar[r]^{a\wedge a}&S^{2m}\ar@/^20pt/[rr]^{1_{S^{2m}}}_{\;}="c"\ar@/_20pt/[rr]_{1_{S^{2m}}}^{\;}="d"
&&S^{2m}\ar@{=>}"a";"b"^{\omega^{\frac{n}{2}}}\ar@{=>}"c";"d"^{\omega^{\frac{m}{2}}}}$$
therefore $2(\alpha\cup_1\alpha)$ corresponds to the self-track
\begin{equation*}\tag{a}
((\omega^{\frac{m}{2}})(a\wedge a))\vc((a\wedge a)(\omega^{\frac{n}{2}})).
\end{equation*}
The self-track $\omega^{\frac{m}{2}}(a\wedge a)$ corresponds to the homotopy class
\begin{equation*}\tag{b}
(\frac{m}{2}(\S^{2(m-1)}\eta))(\S(\alpha\wedge\alpha)).
\end{equation*}
Since 
$\S(\alpha\wedge \alpha)=\pm(\S^{m+1}\alpha)(\S^{n+1}\alpha)$ which is a composite of two triple suspensions
(b) is 
\begin{equation*}\tag{c}
(\alpha\wedge\alpha)(\frac{m}{2}(\S^{2(n-1)}\eta)).
\end{equation*}
Moreover, the self-track $(a\wedge a)(\omega^{\frac{n}{2}})$ corresponds to 
\begin{equation*}\tag{d}
(\alpha\wedge\alpha)(\frac{n}{2}(\S^{2(n-1)}\eta)),
\end{equation*}
so the self-track (a) corresponds to the sum of (c) and (d)
\begin{equation*}
(\alpha\wedge\alpha)(\frac{n+m}{2}(\S^{2(n-1)}\eta)).
\end{equation*}
\end{proof}

\bibliographystyle{amsalpha}
\bibliography{Fernando}

\providecommand{\bysame}{\leavevmode\hbox to3em{\hrulefill}\thinspace}
\providecommand{\MR}{\relax\ifhmode\unskip\space\fi MR }
\providecommand{\MRhref}[2]{%
  \href{http://www.ams.org/mathscinet-getitem?mr=#1}{#2}
}
\providecommand{\href}[2]{#2}
\begin{thebibliography}{BM05b}

\bibitem[Bau91]{ch4c}
H.-J. Baues, \emph{Combinatorial {H}omotopy and 4-{D}imensional {C}omplexes},
  Walter de Gruyter, Berlin, 1991.

\bibitem[BB58]{hcefm}
W.~D. Barcus and M.~G. Barratt, \emph{On the homotopy classification of the
  extensions of a fixed map}, Trans. Amer. Math. Soc. \textbf{88} (1958),
  57--74.

\bibitem[BJ01]{slort}
H.-J. Baues and M.~Jibladze, \emph{Suspension and loop objects and
  representability of tracks}, Georgian Math. J. \textbf{8} (2001), no.~4,
  683--696.

\bibitem[BJM83]{K1p}
M.~G. Barratt, J.~D.~S. Jones, and M.~E. Mahowald, \emph{The {K}ervaire
  invariant problem}, Proceedings of the Northwestern Homotopy Theory
  Conference (Evanston, Ill., 1982) (Providence, RI), Contemp. Math., vol.~19,
  Amer. Math. Soc., 1983, pp.~9--22.

\bibitem[BJP05]{qaI}
H.-J. Baues, M.~Jibladze, and T.~Pirashvili, \emph{Quadratic algebra of square
  groups}, Preprint of the Max-Planck-Institut f\"ur Mathematik MPIM2006-9,
  \texttt{http://arxiv.org/abs/math.CT/0601777}, 2005.

\bibitem[BM05a]{2hg1}
H.-J. Baues and F.~Muro, \emph{Secondary homotopy groups}, Preprint of the
  Max-Planck-Institut f\"ur Mathematik MPIM2006-36,
  \texttt{http://arxiv.org/abs/math.AT/0604029}, 2005.

\bibitem[BM05b]{2hg3}
\bysame, \emph{Smash products for secondary homotopy groups}, Preprint of the
  Max-Planck-Institut f\"ur Mathematik MPIM2006-38,
  \texttt{http://arxiv.org/abs/math.AT/0604031}, 2005.

\bibitem[BM06]{2hg4}
\bysame, \emph{Secondary algebras associated to ring spectra}, In preparation,
  2006.

\bibitem[BP99]{ecg}
H.-J. Baues and T.~Pirashvili, \emph{Quadratic endofunctors of the category of
  groups}, Advances in Math. \textbf{141} (1999), 167--206.

\bibitem[BtD85]{rclg}
T.~Br{\"o}cker and T.~tom Dieck, \emph{Representations of compact {L}ie
  groups}, Graduate Texts in Mathematics, vol.~98, Springer-Verlag, New York,
  1985.

\bibitem[Con84]{2cm}
D.~Conduch\'e, \emph{Modules crois\'es g\'en\'eralis\'es de longueur $2$}, J.
  Pure Appl. Algebra \textbf{34} (1984), no.~2-3, 155--178.

\bibitem[Ell93]{csch}
G.~Ellis, \emph{Crossed squares and combinatorial homotopy}, Math. Z.
  \textbf{214} (1993), no.~1, 93--110.

\bibitem[HM93]{c1Tb}
K.~A. Hardie and M.~E. Mahowald, \emph{On cup-one and {T}oda bracket},
  Quaestiones Math. \textbf{16} (1993), no.~2, 193--207.

\bibitem[Nor90]{norrie}
K.~Norrie, \emph{Actions and automorphisms of crossed modules}, Bull. Soc.
  Math. France \textbf{118} (1990), no.~2, 129--146.

\bibitem[Rut67]{hcmifs}
John~W. Rutter, \emph{A homotopy classification of maps into an induced fibre
  space}, Topology \textbf{6} (1967), 379--403.

\bibitem[Sch11]{dsag}
I.~Schur, \emph{\"{U}ber die {D}arstellung der symmetrischen und der
  alternierenden {G}ruppe durch gebrochene lineare {S}ubstitutionen}, J. Reine
  Angew. Math. \textbf{139} (1911), 155--250.

\bibitem[Ser84]{iwf}
J.-P. Serre, \emph{L'invariant de {W}itt de la forme {${\rm Tr}(x\sp 2)$}},
  Comment. Math. Helv. \textbf{59} (1984), no.~4, 651--676.

\bibitem[Tei92]{t4mffg}
P.~Teichner, \emph{Topological four-manifolds with finite fundamental group},
  Dissertation zur {E}rlangung des {G}rades {D}oktor der {N}aturwissenschaften,
  Johannes-Gutenberg Universit\"at in Mainz, Mainz, M\"arz 1992.

\bibitem[Tod62]{toda}
H.~Toda, \emph{Composition methods in homotopy groups of spheres}, Annals of
  Mathematics Studies, No. 49, Princeton University Press, Princeton, N.J.,
  1962.

\end{thebibliography}
\end{document}